\documentclass[12pt]{amsart}

\usepackage{times,amsfonts,amsmath,amstext,amsbsy,amssymb,
  amsopn,amsthm,upref,eucal, amscd} 
\usepackage[T1]{fontenc}

\newtheorem{theorem}{Theorem}[section]
\newtheorem{lemma}[theorem]{Lemma}
\newtheorem{corollary}[theorem]{Corollary}

\numberwithin{equation}{section}

\theoremstyle{definition}

\newtheorem{remark}[theorem]{Remark}

\newcommand{\field}[1]{\mathbb{#1}}
\newcommand{\R}{\field{R}} 
\newcommand{\N}{\field{N}}
\newcommand{\C}{\field{C}} 
\newcommand{\Z}{\field{Z}}
\newcommand{\Q}{\field{Q}} 
\newcommand{\h}{\field{H}}
\newcommand{\Sone}{\field{S}^1}
\newcommand{\Ab}{{\hat \Omega}_g}
\newcommand{\WAb}{\hat W}
\newcommand{\Teich}{\tilde{\Omega}_g}
\newcommand{\HTeich}{\tilde{H}}
\newcommand{\WTeich}{\tilde{W}}
\newcommand{\dc}{\hat\Phi}
\newcommand{\Cal}{\mathcal}

\newcommand{\<}{{\langle}} 
\renewcommand{\>}{{\rangle}}
\newcommand{\pref}[1]{(\ref{#1})}

\begin{document}
\title{Deviation of ergodic averages for rational polygonal billiards}
\author{J.~S.~Athreya and G.~Forni}
\email{jathreya@math.princeton.edu}
\email{gforni@math.umd.edu}
\address{Dept. of Mathematics, Princeton University, Princeton, NJ}
\address{Dept. of Mathematics, University of Maryland, College Park, MD}
\thanks{J.S.A. supported by NSF grant DMS
    0603636.}
\thanks{G.F. supported by NSF Grant DMS 0244463.}    
\begin{abstract}
We prove a polynomial upper bound on the deviation of ergodic averages for almost all directional 
flows on every translation surface, in particular, for the generic directional flow of billiards in any Euclidean polygon with rational angles.
\end{abstract}
\maketitle

\section{Introduction and Main Results}\label{intro}
In a celebrated paper \cite{KMS}, S. Kerckhoff, H. Masur and J. Smillie proved the unique ergodicity
for the directional flow of a rational polygonal billiard in almost all directions.  In this paper we prove a power-law upper bound on the speed of ergodicity for weakly differentiable functions. 

Similar bounds for hyperbolic systems are usually modeled on the \emph{Law of Iterated Logarithms }
(or on the Central Limit Theorem) for independent stochastic processes. In such cases the deviation 
exponent for \emph{ergodic integrals }is  universal, equal to $1/2$. Among well-known models  of non-hyperbolic behaviour, bounds on the speed of ergodicity are provided by the \emph{Denjoy-Koksma inequality }for rotations of the circle or, equivalently, for linear flows on the $2$-torus. In this case, for generic systems the deviation of ergodic integrals is at most logarithmic in time.

There are a few results known for dynamical systems with  \emph{intermediate} behaviour, that is, systems which display slow (polynomial) divergence of nearby orbits. Such systems have been called
\emph{parabolic} (see for instance \cite{KHHB}). An important example is the horocycle flow on the unit tangent bundle of a surface of constant negative curvature. M. Burger \cite{Burger} proved polynomial upper bounds for geometrically finite surfaces and found that the deviation exponent depends on the spectrum of the Laplace-Beltrami operator. His results were later strengthened in \cite{FF1} where a quite refined asymptotics, including lower bounds, was proved in the case of compact surfaces. Another class of examples where polynomial upper bounds have been proved is given by nilflows on Heisenberg manifolds \cite{FF2}. In all the above parabolic examples, the generic system is \emph{uniquely ergodic} and the polynomial upper bounds are \emph{uniform} with respect to the initial conditions, in sharp contrast to the hyperbolic situation.

The dynamics of a rational polygonal billiard can be viewed as the geodesic flow on a flat surface
with trivial holonomy, called a translation surface. In fact, translation surfaces arising from billiards
are rather special among higher genus translation surfaces. Flat $2$-dimensional tori are the only translation surfaces of genus $1$ and correspond to \emph{integrable }billiard tables, such as rectangles or the equilateral triangle. In this case orbits do not diverge and the system is called
\emph{elliptic}. In the higher genus case, the geodesic flow is \emph{pseudo-integrable} in the sense that the phase space is foliated by surfaces of higher genus. Each invariant surface is determined by fixing an angle of the unit tangent vector with respect to the horizontal direction.  Since the surface has a translation structure (the flat metric has trivial holonomy) the horizontal direction is well-defined and the angle is invariant under the geodesic flow. The divergence of nearby orbits is
produced by the presence of (conical) singularities which split bundles of nearby trajectories. 
By the Kerckhoff-Masur-Smillie theorem, the flow is uniquely ergodic on almost all invariant surfaces. In this paper we prove a polynomial bound on the speed of ergodicity for a (smaller) full measure set of invariant surfaces. 

The flow on each invariant surface, called the directional flow, is closely related to \emph{interval
exchange transformations} (IET's). In fact, all Poincar\'e return maps on transverse intervals are
IET's. For a generic IET a polynomial bound on the deviation of ergodic sums for certain piece-wise
constant functions was proved by Zorich in his groundbreaking work \cite{Zorich1}, \cite{Zorich2}, \cite{Zorich3}, \cite{Zorich4}.
In terms of translation flows, Zorich's result yields a power-law upper bound for the deviation in
homology of the generic directional flow on the generic surface (see Theorem \ref{mainhom} below).
A similar power-law upper bound on the deviation of ergodic integrals of \emph{weakly differentiable functions}, together with results on lower order deviations and lower bounds on deviations, were
proved by the second author in \cite{Forni}, after conjectures of Zorich \cite{Zorich1} (see also \cite{Zorich4}) and Kontsevich \cite{K}. We emphasize that in genus $g\geq 2$ all of the above mentioned results on deviations are for \emph{generic }translation surfaces and \emph{tell us absolutely nothing on rational billiards}, which form a set of measure zero. This situation is quite unsatisfactory, since the study of billiards and of related mechanical systems is one of the main motivations of the theory. Our current result is a first step in addressing this gap.

\medskip
\subsection{Holomorphic Differentials} Let $M$ be a Riemann surface of genus $g\geq 1$. 
Let $\text{Hol}(M)$ denote the space of holomorphic differentials on $M$, that is, the space
of all tensors of the form $f(z)dz$ in local coordinates. Any holomorphic differential $\omega \in \text{Hol}(M)$ determines a unique flat metric with conical singularities at the zeros of the holomorphic differential. Given $\omega \in \text{Hol}(M)$, one obtains (via integration of the form) an atlas of charts to $\C \cong \R^2$, with transition maps of the form $z \mapsto  z + c$. Viceversa, given such an atlas of charts, one obtains a holomorphic differential by pulling back the form $dz$ on $\C$. For this reason pairs $(M,\omega)$ are also called \emph{translation surfaces} (this terminology was first introduced 
by E. Gutkin and C. Judge in \cite{GuJu}).  Any differential $\omega\in \text{Hol}(M)$ also determines a pair of transverse oriented measured foliations, defined by $\{\operatorname{Re} (\omega) =0\}$ \emph{(vertical foliation)}, $\{\operatorname{Im} (\omega) = 0\}$ \emph{(horizontal foliation)}. These foliations have saddle-like singularities (possibly degenerate) at the zero of the holomorphic differential. The \emph{vertical} and \emph{horizontal flows }are the flows moving with unit speed along the leaves of the vertical and horizontal foliations respectively, in the positive direction. More in general, a \emph{directional flow }is the flow moving with unit speed in a direction forming a fixed angle with the (positive) vertical direction. We are interested in studying the geodesic flow of the flat metric associated with any $\omega\in \text{Hol}(M)$. It is well-known that the unit cotangent bundle of $M$ is foliated with invariant surfaces of genus $g\geq 1$ and that the restriction of the geodesic flow to an invariant surface is isomorphic to a directional flow on $M$.  In fact, it is sufficient to study the vertical (or the horizontal) flow for the one-paramer families of holomorphic differentials on $M$ obtained by `rotations' of $\omega \in \text{Hol}(M)$ (see below).

\subsection{Moduli space} 
Let $S= S_g$ be a compact surface of genus $g \geq 2$.  Let $\Omega_g$ be the moduli space of 
unit-area holomorphic differentials on $S$. That is, a point in $\Omega_g$ is an equivalence class of pairs $(M, \omega)$, where $M$ is a Riemann surface of genus $g\geq 1$ (that is, a complex structure on $S$) and $\omega$ is a holomorphic differential on $M$ normalized so that the associated metric has unit total area.  The equivalence relation is defined as follows: two pairs $(M_1, \omega_1)$ and $(M_2, \omega_2)$ are equivalent if there is a biholomorphic map $f:M_1 \rightarrow M_2$ such that $f^{*} \omega_1 = \omega_2$. For notational convenience, we will write points in $\Omega_g$ simply as $\omega$ and $M_\omega$ will denote the underlying Riemann surface.

The space $\Omega_g$ is an orbifold vector bundle over $\Cal R_g$, the moduli space of Riemann surfaces. 
The fiber over each point $M \in \Cal R_g$ is the vector space of holomorphic differentials on $M$. There is a natural stratification of $\Omega_g$ by integer partitions of $2g-2$: each stratum can be described as the subset of $\Omega_g$ formed by all differentials with a given pattern of zero multiplicities. Strata are never compact and not always connected. However, they have at most finitely many connected components~\cite{KZ}. In our paper we will work with one of these connected components, call it $\Omega$.

There is a naturally defined action of the group $SL(2, \R)$  of $2\times 2$ matrices with determinant
$1$ on the moduli space $\Omega_g$ of holomorphic differentials. In fact, $SL(2,\R)$ acts linearly
on $\R^2 \equiv \C$ and its action on $\Omega_g$ can be defined by post-composition 
with holomorphic differentials or, equivalently, with charts in any atlas for the corresponding translation structure. Given any translation surface $(M,\omega)$, post-composition with any $A\in SL(2, \R)$ uniquely defines a translation surface $(M_A, \omega_A)$. The $SL(2,\R)$ action preserves the multiplicities of zeroes of the holomorphic differentials (or equivalently the angles of the conical singularities), hence it preserves each stratum of the moduli space.

We are interested in the dynamics of the vertical flow for almost all holomorphic differentials
in every orbit of the standard maximal compact subgroup $SO(2,\R)$,

$$SO(2,\R)=\left\{r_{\theta} = \left( \begin{array}{cc} \cos \theta & \sin \theta \\
-\sin \theta & \cos \theta \end{array}\right):  \theta \in \Sone \right\} \,.$$ 

The action of $SO(2,\R)$, known as the \emph{circle flow}, preserves the underlying holomorphic structure (as well as the flat metric), so it acts as the identity when projected to $\Cal R_g$. From a geometrical point of view, the matrix $r_{\theta}$ simply rotates the vertical direction by angle 
$\theta\in \Sone$.

\subsection{Main Theorem}

Before we state our main result, we fix some more notation. If we fix $\omega \in \Omega$, let $\varphi_t$ denote the vertical flow. Let $\varphi_{\theta, t}$ denote the vertical flow associated to $r_{\theta} \omega$ (this is simply the directional flow in the direction at angle $\theta\in \Sone$ with the positive vertical direction of $\omega$). We say a point $x \in S$ is non-singular for $\theta\in \Sone$ if it is not on a singular leaf of the vertical foliation associated to $r_{\theta} \omega$. Let $A_{\omega}$ denote the area form on $S$ associated to $\omega$ (note that $A_{r_{\theta}\omega} = A_{\omega}$). 

\begin{theorem}\label{main}  There is an $\alpha = \alpha (\Omega) >0$  such that the following holds. For all $\omega \in \Omega$ there is a measurable function $K_\omega: \Sone \to \R^+$ such that
 for almost all $\theta \in \Sone$ (with respect to Lebesgue measure), for all functions $f $ in the standard Sobolev space $H^1(S)$ and for all non-singular $x \in S$,  
\begin{equation}
\label{mainresult} 
\left \vert \int_{0}^{T} f(\varphi_{\theta, t}(x)) dt  - T\int f dA_\omega \right \vert   \,\leq \, K_\omega(\theta)
\Vert f \Vert_{H^1(S)} \, T^{1-\alpha}  \,.
\end{equation}
\end{theorem}

We also have the following result for the growth of homology classes: Fix $\omega\in \Omega$, and let $x\in S$ be non-singular. For $T>0$, define $h_{\omega, x}(T) :=[ \bar{\gamma}_{\omega, x} (T) ]\in H_1(S,\R)$, where $\bar{\gamma}_{\omega, x}(T)$ is the closed curve given by taking the piece of leaf $\gamma_{\omega, x}(T) := \{\varphi_t(x)\}_{0}^{T}$, and `closing' it up by connecting the endpoints by 
any given curve of bounded length. We remark that $h_{\omega,x}(T)$ depends on the choice
of the curve joining the endpoints, but it is uniquely defined up to
addition of a uniformly bounded term.

\begin{theorem}\label{mainhom} There is an $\alpha = \alpha (\Omega) >0$ (the same as in Theorem~\ref{main}) and for all $\omega \in \Omega$ there is a measurable function $K'_\omega: \Sone \to \R^+$, such that for almost all $\theta \in \Sone$, there is a homology class $h_{\theta} \in H_1(S,\R)$ so that  for all non-singular $x \in S$, 
\begin{equation}
\label{homresult}
 \left \vert h_{r_{\theta}\omega, x}(T) - h_{\theta} T \right \vert \,  \leq \, K'_\omega(\theta) \, T^{1-\alpha}\,.
 \end{equation}
 \end{theorem}
 
 Interestingly, while generically the deviation exponent in Theorem \ref{main} or Theorem 
 \ref{mainhom} is positive \cite{Forni}, there are higher genus surfaces for which it can be zero, and in~\cite{ForniSurvey}, an example is given of a surface of genus $3$ where the exponent is zero in almost all directions. Thus, we cannot prove a general non-trivial \emph{lower} bound  for the deviation 
 exponent. Unpublished work of the second author suggests that such non-generic examples do not exist in genus $2$. More precisely, it can be proved that in genus $2$ the (upper) second Lyapunov exponent of the so-called Kontsevich-Zorich cocycle (see \S \ref{KZcocycle} and \S \ref{KZgap}) is strictly positive for all holomorphic differentials and for almost all angles. In particular, the second Lyapunov exponent is strictly positive with respect to all $SL(2,\R)$-invariant measures.  Recently, M. Bainbridge \cite{Bainbridge} has derived explicit numerical values for the (second) exponent, with respect to all $SL(2,\R)$-invariant measures, from a formula found by M. Kontsevich  \cite{K} (see also \cite{Forni} for a complete proof of the formula).  It turns out that the values are equal to $1/3$ for all measures supported in the stratum corresponding to a double zero of the holomorphic differential and to $1/2$ for all measures in the stratum corresponding to two simple zeros. These values were conjectured by Kontsevich and Zorich (see for instance \cite{K} for the case of the absolutely continuous invariant measures on each stratum) but their calculations have remained unpublished. 
    
 \subsection{Rational Polygonal Billiards}
 \label{RPB}
A special class of holomorphic differentials (or translation surfaces) is given by \emph{rational polygonal billiards}. Let $P \subset \R^2$ be a Euclidean polygon and let $G(P)\subset O(2,\R)$ be the subgroup generated by all reflections with axis parallel to an edge of $P$ (and passing through the origin). The polygon $P$ is called \emph{rational} if $G(P)$ is finite. A necessary condition, which is also sufficient if $P$ is simply connected, is that the angles of $P$ belong to $\pi \Q$ (see, for example, the excellent survey~\cite{MasurTab}, \S 1.3). The billiard flow on $P$ is a discontinuous Lagrangian (Hamiltonian) flow on the unit tangent bundle $T_1(P)\equiv P\times \Sone$. The trajectory of any $(x,v) \in T_1(P)$ moves with unit speed along a straight line in the direction $v\in \Sone$ up to the boundary $\partial P$ where it is reflected according to the law of geometric optics (\emph{angle of incidence equal angle of reflection}), which follows from the assumption that collisions with the boundary are elastic. 

The billiard flow on a rational table $P$ leaves invariant the angle function $\Theta: P\times \Sone \to 
\R$ obtained as a composition of the canonical projections 
\begin{equation*}
\Theta: P\times \Sone \rightarrow    \Sone \rightarrow    \Sone/G(P) 
\end{equation*}
onto the quotient $ \Sone/G(P)$, which can be identified to a compact interval $I(P) \subset \R$.
It follows that the phase space $P\times \Sone$ is foliated by the level surfaces $S_{P,\theta} = \{ (x,v) \in P\times \Sone \vert \Theta(x,v) =\theta\}$ which are invariant under the billiard flow and have natural translation structures induced by the translation structure on $P \subset \R^2$. By following the unfolding construction of Zemljakov-Katok~\cite{KZ1}, it is possible to show that all invariant translation surfaces $S_{P,\theta}$ for $\theta
\not \in \Theta^{-1} (\partial I(P))$ can be identified with rotations of a fixed translation surface $S_P$ of genus $g(P) \geq 1$ so that the billiard flow restricted to $S_{P,\theta}$ can be identified with the directional flow on $S_P$ in the direction at angle $\theta\in \Sone$ from the vertical (see the original work by E. Gutkin \cite{Gu} or the survey by H. Masur and S. Tabachnikov \cite{MasurTab}, \S 1.5). 

If $P$ is rational billiard table, let $A_{P, \theta}$ denote the area form on the invariant translation surface $S_{P,\theta}$ for all $\theta \in \Sone$. We have the following corollary:

\begin{corollary}  
\label{ratbill}
Let $P\subset \R^2$ be a rational polygon. For any $\theta \in \Sone$, let 
$\psi_{t, \theta}$ be the restriction of the billiard flow to the invariant surface $S_{P,\theta} \subset P \times \Sone$.  There exist an $\alpha = \alpha(P)>0$  (depending only on the shape of $P$, in particular on the stratum arising from the unfolding procedure) and a measurable function $K_P:\Sone\to \R^+$ such that, for almost all $\theta\in \Sone$, for all $f$ in the standard Sobolev space $H^1(S_{P,\theta})$, and all $x\in S_{P,\theta}$ for which $\psi_{\theta, t}(x)$ is defined for all $t >0$,
\begin{equation}
\left | \int_{0}^{T} f(\psi_{\theta, t}(x)) dt  - T\int f dA_{P,\theta} \right | \,\leq \, K_P(\theta) \Vert f \Vert_{H^1(S_{P,\theta}) } \, T^{1-\alpha} \,.
\end{equation}
\end{corollary}

The first estimates on the speed of ergodicity for billiards in polygons were obtained 
by Ya. B. Vorobets in \cite{Vorobets1}, \cite{Vorobets2}. For rational polygonal billiards his bound
 (which holds for Lipschitz functions vanishing on the boundary of the billiard table) is much weaker 
 than ours. In fact, it is far from being polynomial. However, he can control explicitly the mean over the angle of the deviation of ergodic averages in terms of the shape of the billiard table, hence his result yields ergodicity for an explicit class (measure zero, but topologically large) of non-rational polygonal billiards by the approximation method of Katok and Zemljakov \cite{KZ1} (see also \cite{KMS}). 
 Vorobets methods are not based on Teichm\"uller theory.  We remark that the dependence of the exponent in our Corollary \ref{ratbill} with respect to the shape of the billiard table is not sufficiently  explicit to derive by approximation an effective ergodicity result for non-rational polygonal billiards.

\medskip

The rest of this paper is structured as follows: In section~\ref{dict}, we discuss the Teichm\"uller geodesic flow and the dictionary between ergodic properties of foliations (or flows) on surfaces and recurrence properties of Teichmu\"uller orbits, as well as recalling the results from~\cite{Forni} for deviations of \emph{generic} vertical flows. In section~\ref{cocycles}, we discuss the Kontsevich-Zorich cocycle, a symplectic cocycle over the flow, and distributional generalizations introduced in \cite{Forni}. In section~\ref{spectral}, we prove the key lemma for our theorems, an estimate on the spectral gap 
of the Kontsevich-Zorich cocycle and of the distributional cocycle, which combines explicit formulas from~\cite{Forni} with large deviations estimates on the Teichm\"uller flow from~\cite{Athreya}. Finally, in section~\ref{proofs}, we bring these results together in order to prove Theorems~\ref{main} and~\ref{mainhom}.

\section{A Dictionary}\label{dict}

 Let 
 \begin{equation}
 \label{eq:Tflow}
 A = \left\{g_t = \left(\begin{array}{cc} e^{t} & 0 \\ 0 & e^{-t}
\end{array}\right): t \in \R\right\}\,.
\end{equation}

The action of $A$ on a stratum $\Omega$ is known as the \emph{Teichm\"uller geodesic flow},
since the projection of any $A$-orbit yields a geodesic in the Teichm\"uller metric on $\Cal R_g$ (and in fact, all Teichm\"uller geodesics arise this way). In terms of foliations, if we write $\omega = (\operatorname{Re}(\omega), \operatorname{Im}(\omega))$, we have 
$$
g_t \omega = (e^{t}\operatorname{Re}(\omega), e^{-t}\operatorname{Im}(\omega))\,.
$$
Each stratum $\Omega$, while non-compact, is endowed with a canonical absolutely continuous, ergodic, $A$-invariant (in fact, $SL(2,\R)$-invariant) probability measure $\mu = \mu_{\Omega}$
\cite{Masur1}. In any stratum $\Omega$ of surfaces of genus $g\geq 2$ the set of points arising from rational billiards has $\mu$-measure $0$. As a consequence, results for $\mu$-generic points which are easier to obtain by methods of ergodic theory, do not directly apply to billiards. As we have remarked above, this is a serious difficulty in the study of rational polygonal billiards.

There is a well-studied dictionary between the ergodic properties of the vertical flow (or foliation) associated to $\omega \in \Omega$, and the recurrence properties of the forward geodesic trajectory $\{g_t \omega\}_{t \geq  0}$ (a similar discussion can be had about the horizontal flow and the backward trajectory $\{g_t \omega\}_{t \le 0}$). The first main result in this dictionary is known as \emph{Masur's Criterion}:

\begin{theorem}(Masur~\cite{Masur2}) If the vertical foliation $\{\operatorname{Re}(\omega) = 0\}$ is non-uniquely ergodic (that is, there is more than one transverse invariant probability measure), $\{g_t \omega\}_{t\geq 0}$ is divergent in $\Omega$.\end{theorem}

Combining this with the fact that $\{g_t\}$ is ergodic with respect to $\mu$~\cite{Masur1, Veech1}, one obtains that for $\mu$-almost every $\omega \in \Omega$, the vertical foliation is uniquely ergodic. However, since (the set of surfaces arising from) billiards have measure zero, this is not useful for billiards.

This difficulty was resolved by Kerckhoff, Masur and Smillie~\cite{KMS}, who analyzed the recurrence behavior of $\{g_t r_{\theta} \omega\}_{t \geq 0}$ for a fixed $\omega$ and $\theta\in \Sone$ varying. The main result in~\cite{KMS} is as follows:

\begin{theorem}(Kerckhoff, Masur and Smillie~\cite{KMS}) For any $\omega \in \Omega$, the set of $\theta\in \Sone$ for which $\{g_t r_{\theta} \omega\}_{t \geq 0}$ is divergent has measure $0$.\end{theorem}

Thus, by Masur's criterion, for all translation surfaces (and in particular for all rational polygonal billiards), and in almost all directions, the vertical flow is uniquely ergodic. 

As mentioned above, the first results on the \emph{deviation} of ergodic averages for such systems were proved by Zorich~\cite{Zorich3}, who proved a result equivalent to Theorem~\ref{mainhom} for \emph{generic} interval exchange transformations. Zorich \cite{Zorich1}, \cite{Zorich4} and Kontsevich \cite{K} conjectured precise power-laws for the deviations of ergodic averages of smooth functions for \emph{generic }interval exchange transformations and translations flows. The Kontsevich-Zorich conjectures, including Theorem~\ref{main} in the \emph{generic} case, were proved by the second author~\cite{Forni}, with exception of the \emph{simplicity }of the deviation (Lyapunov) spectrum, recently proved by A. Avila and M. Viana \cite{AV}. Once again, none of these results apply to rational polygonal billiards. 

The results of~\cite{Forni} were obtained by careful study of the Lyapunov exponents of the Kontsevich-Zorich cocycle and of its distributional extensions over the Teichm\"uller flow, which we discuss in the next section. In particular, these results are based on the recurrence behavior of \emph{generic} Teichm\"uller trajectories to compact sets given by ergodicity. In ~\cite{Athreya} the first author, building on work of Eskin-Masur~\cite{EskinMasur}, analyzed
in depth the recurrence to a certain exhaustion by compact subsets of the moduli space for almost 
all geodesic trajectories in \emph{every }orbit of the group $SO(2,\R)$.  We will use the following key result:

\begin{theorem}\label{proportion} For any $\eta>0$, there is a compact set $C = C(\eta) 
\subset \Omega$ such that for all $\omega \in \Omega$ and almost all $\theta\in \Sone$,
\begin{displaymath}\limsup_{T \rightarrow \infty} \frac{1}{T} \left|\{ 0 \le t \le T: g_t
r_\theta \omega \notin C\}\right| \le
\eta\,.\end{displaymath}
\end{theorem}

\vspace{.1in}

\textbf{Remark:} This is essentially Corollary 2.4 from~\cite{Athreya}. However, there the compact set $C$ depends on the basepoint $\omega\in \Omega_g$. We thank Barak Weiss for pointing out that the proof of this result in~\cite{Athreya} shows that for any stratum $\Omega
\subset \Omega_g$, the set $C$ can be picked independently of $\omega \in \Omega$.

\section{Cocycles over the Teichm\"uller flow}
\label{cocycles}

\subsection{The Kontsevich-Zorich cocycle}
\label{KZcocycle}

 The Kontsevich-Zorich cocycle is a multiplicative symplectic cocycle over the Teichm\"uller
geodesic flow on the moduli space of holomorphic abelian differentials on compact Riemann surfaces. This cocycle appears in the study of the dynamics of interval exchange transformations and of translation flows on surfaces,  for which it represents a renormalization dynamics, and of the 
Teichm\"uller flow itself. In fact, the study of the tangent cocycle of the Teichm\"uller flow can be reduced 
to that of the Kontsevich-Zorich cocycle.

 Let $\Teich$ be the Teichm\"uller space of abelian differentials on a closed surface 
$S$ of genus $g\geq 1$. Since points in $\Teich$ are pairs $(M, \omega)$, with $M$ a marked Riemann surface, the trivial cohomology bundle $\Teich \times H^1(S,\R)$ is well-defined. The {\it Kontsevich-Zorich cocycle }$\{\Phi_t\}_{t\in\R}$ can be defined as the projection of the trivial cocycle
\begin{equation}
\label{eq:KZcocycle}
g_t \times \hbox{id}: \Teich\times H^1(S,{\R})\to \Teich\times H^1(S,{\R}) 
\end{equation}
onto the orbifold vector bundle ${\Cal H}^1_g(S,{\R})$ over $\Omega_g$ defined as 
\begin{equation}
\label{eq:CB}
{\Cal H}^1_g(S,{\R}):=\bigl(\Teich\times H^1(S,{\R})\bigr) /\Gamma_g\,.
\end{equation}
The mapping class group $\Gamma_g$ acts naturally on the bundle $\Teich\times H^1(S,{\R}) $ since it acts naturally on the real cohomology $H^1(S,\R)$ by pull-back.

 The cohomology bundle ${\Cal H}^1_g(S,{\R})$ can be endowed with the structure of a smooth euclidean bundle with respect to the Hodge inner product.  In fact, by the Hodge theory on Riemann surfaces (~\cite{FK92}, III.2), any real cohomology class $c\in H^1(M, \R)$ can be represented as the real (or imaginary) part of a holomorphic differential $h \in \text{\rm Hol}(M)$ on a Riemann surface $M$.  Let $\omega\in \Omega_g$ and let $M_\omega$ be the underlying Riemann surface. For any $c \in H^1(M_\omega, \R)$, its  Hodge norm is defined as 
\begin{equation}
\label{eq:Hodgenorm}
\Vert c \Vert^2_{\omega} := \frac{i}{2}\,  \int_S  h \wedge  \overline {h}  \quad \text{ \rm if }\,\,
c=[\operatorname{Re} (h)]\,, \,\, h\in \text{\rm Hol}(M_\omega)\,.
\end{equation}
 We remark that the Hodge norm is defined in terms of the Riemann surface $M_\omega$ but it is otherwise independent of the differential $\omega\in \Omega_g$.
 
 Real cohomology classes on $S$ can also be represented in terms of meromorphic 
functions in $L^2_{\omega}(S)$ (see \cite{Forni}, \S 2). In fact, any holomorphic differential $\omega\in \Omega_g$ induces an isomorphism between the space $\text{\rm Hol}(M_\omega)$ of all holomorphic differentials on $M_\omega$ and the subspace ${\Cal M}_{\omega}$ of meromorphic functions in $L^2_{\omega}(S)$. Such a subspace can be characterized as the space of all meromorphic functions with poles at the zeros of $\omega$ of orders bounded in terms of the order of the corresponding zero \cite{Forni1}, \cite{Forni}. 

Any holomorphic differential $h\in\text{\rm Hol}(M_\omega)$ can be written in terms of a meromorphic function $m\in L^2_{\omega}(S)$ as follows:
\begin{equation}
\label{eq:holmer}
h :=  m \,\omega \,, \quad m \in \Cal M_{\omega}\,.
\end{equation}
 The following \emph{representation of real cohomology classes} therefore holds:
\begin{equation}
\label{eq:repr}
c\in H^1(S,\R) \Longleftrightarrow  c=[ \operatorname{Re} ( m \, \omega)]  \,, \quad m \in \Cal M_{\omega}.
\end{equation}
The map $c_{\omega}:{\Cal M}_{\omega} \to H^1(S,{\R})$ given by the representation \pref{eq:repr} is bijective and it is in fact {\it isometric }if the space ${\Cal M}_{\omega}$ is endowed with the Euclidean structure induced by $L^2_{\omega}(S)$ and $H^1(S,{\R})$ with the {\it Hodge norm }$\Vert \cdot \Vert_\omega$ introduced in \pref{eq:Hodgenorm}. In fact,  the following identity holds:
\begin{equation}
\label{eq:Hodgenormbis}
\Vert c_{\omega}(m) \Vert^2_\omega :=  \int_S  \vert m \vert ^2 \,dA_{\omega} \, ,  \quad \text{ \rm for all }\,\,
m\in \Cal M_{\omega}.
\end{equation}

The Kontsevich-Zorich cocycle was introduced in \cite{K} as a continuous-time 
version of the Zorich cocycle. The Zorich cocycle was introduced earlier by A. Zorich \cite{Zorich2}, in order to explain polynomial deviations in  the homological asymptotic behavior of 
 typical leaves of orientable measured foliations on compact surfaces, a phenomenon he had discovered in numerical experiments \cite{Zorich1}.  We recall that the real homology $H_1(S,\R)$ 
 and the real cohomology $H^1(S,\R)$ of an orientable closed surface $S$ are (symplectically) isomorphic by Poincar\'e duality.
 
 Zorich proved in \cite{Zorich2} that the integrability condition of Oseledec's multiplicative ergodic theorem is satisfied for a suitable acceleration of the Rauzy-Veech induction, now often called
 the Zorich induction. The integrability condition is immediate for the Kontsevich-Zorich cocycle.
 In fact, the logarithm of the Hodge norm of the cocycle is a bounded function on the moduli 
 space, hence it is integrable with respect to any probability measure.

 Since the Kontsevich-Zorich cocycle is symplectic, its Lyapunov spectrum is symmetric
with respect to the origin. Hence for any $g_t$-invariant ergodic probability measure $\mu$ on $\Omega_g$, the Lyapunov spectrum over $(g_t, \mu)$ is equal to an ordered set of the form:
\begin{equation}
\label{eq:KZexp}
\lambda_1^\mu \geq \dots \geq \lambda^\mu_g \geq \lambda^\mu_{g+1}= -\lambda^\mu_g \geq
\dots \geq \lambda^\mu_{2g} = -\lambda_1^{\mu}\,.
\end{equation}
 A crucial property of the cocycle is the following `spectral gap' theorem, proved by W. Veech \cite{Veech} for a class of measures satisfying certain integrability conditions (including the canonical absolutely continuous $g_t$-invariant measure on any connected component of any stratum of the moduli space) and generalized by the second author \cite{Forni}, \cite{ForniSurvey}.

\begin{theorem} \label{thm:mugap} For any $g_t$-invariant ergodic probability measure $\mu$ on the moduli space $\Omega_g$, the following inequality holds:
\begin{equation}
\label{eq:l2ub}
\lambda^{\mu}_2 <\lambda^{\mu}_1=1 \,\,. 
\end{equation}
\end{theorem}

We recall that Zorich \cite{Zorich1}, \cite{Zorich2}, \cite{Zorich4} and Kontsevich \cite{K} conjectured that, if $\mu$ is he canonical absolutely continuous $g_t$-invariant measure on any connected component of any stratum of the moduli space, then the Lyapunov exponents in \pref{eq:KZexp} are $(a)$ all non-zero (non-uniform hyperbolicity); $(b)$ all distinct (simplicity). Part $(a)$ of the Kontsevich-Zorich
conjecture was first proved by the second author in \cite{Forni}, while part $(b)$, which is stronger, was
proved by a completely different approach by Avila and Viana \cite{AV}. This paper does not rely
on any of the above mentioned results on the Lyapunov spectrum of the Kontsevich-Zorich cocycle.
In fact, our results are based on a stronger version of the spectral gap, Theorem \ref{thm:mugap}, which
will be proved in section  \S \ref{spectral}.

 \subsection{ Weighted Sobolev spaces}
 \label{WSS}
 We recall the definition, introduced in \cite{Forni1} (see also \cite{Forni}, \cite{ForniSurvey}) of the weighted Sobolev spaces associated with an abelian differential.
 
Let $\Ab$ denote the space of all abelian differentials on a closed surface $S$ of genus $g\geq 1$, that is, the space of all complex-valued closed $1$-forms which are holomorphic with respect to some complex structure on $S$.
 
 For any $\omega\in \Ab$, we consider the Sobolev spaces $H^1_\omega(S) \subset L^2_\omega(S)$ naturally associated with the flat metric $R_\omega$ (with conical singularities) induced by $\omega$ on $S$. Such spaces are defined as follows. Let $A_\omega$ be the area form associated with $\omega\in  \Ab$, that is,
 \begin{equation}
 \label{eq:area}
  A_\omega :=    \frac{i}{2} \, \omega \wedge \bar \omega \,.
 \end{equation}
 Let $(X_\omega, Y_\omega)$ be a $R_\omega$-orthonormal system of parallel vector fields
 on $S \setminus \Sigma_\omega$ such that $X_\omega$ is tangent to the horizontal foliation in the positive direction and  $Y_\omega$ is tangent to the vertical foliation, also in the positive direction,
 in other terms $(X_\omega,Y_\omega)$ are uniquely determined by the conditions:
  $$
  \imath_{X_\omega} \omega = 1  \quad \text{ \rm and }  \quad   \imath_{Y_\omega} \omega = i \,.
  $$

 We recall that,  with respect to any local canonical coordinate $z=x+iy$ on $S \setminus \Sigma_\omega$ (such that $\omega=dz$), the metric $R_\omega$ has the form 
 $R_\omega= dx^2 + dy^2$, hence
 \begin{equation}
A_\omega = dx \wedge dy \, ,\quad (X_\omega, Y_\omega) = (\frac{\partial}{\partial x},\frac{\partial}{\partial y}) \,,
\end{equation}
 while with respect to a canonical coordinate centered at a zero $p\in S$ of $\omega$ of order $n$ 
 (such that $\omega =z^n dz$), it has the form $R_\omega = \vert z\vert^{2n} (dx^2 + dy^2)$, hence
    \begin{equation}  
  \begin{aligned}
  A_\omega&= \vert z\vert^{2n} dx \wedge dy   \,,\\
 X_\omega &=  \frac{1}{\vert z\vert^{2n}}\left( \operatorname{Re}(z^n)\, \frac{\partial}{\partial x} - \operatorname{Im}(z^n)\, \frac{\partial}{\partial y}\right)\,, \\
 Y_\omega &=   \frac{1}{\vert z\vert^{2n}}\left(\operatorname{Im}(z^n) \, \frac{\partial}{\partial x} + \operatorname{Re}(z^n)\, \frac{\partial}{\partial y} \right)\,.
 \end{aligned}
\end{equation}

The weighed Sobolev spaces $H^1_\omega(S)\subset L^2_\omega(S)$ are then
defined as follows: the space $L^2_\omega (S) := L^2(S, dA_\omega)$ is the standard space
of square-integrable functions with respect to the area and 
\begin{equation}
H^1_\omega (S) := \{ f \in L^2_\omega(S) \vert   X_\omega f\,  \text{ \rm and }\, Y_\omega f \in 
L^2_\omega (S)\}\,.
\end{equation}
The space $L^2_\omega(S)$  and its norm, defined for all $f \in L^2_\omega(S)$ as
\begin{equation}
\vert f \vert_{0,\omega} :=  \left ( \int_S  \vert f \vert^2 \, dA_\omega \right)^{1/2} \,,
\end{equation}
do depend on the abelian differential $\omega \in \Ab$. However,
since the area form $A_\omega$ on $S$ is invariant under the natural action of $SL(2, \R)$ on
$\Ab$, the space $L^2_\omega(S)$  and its norm are also $SL(2,\R)$-invariant
(in particular, $g_t$-invariant).

 It was proved in \cite{Forni1}, \cite{Forni} that, while $L^2_\omega(S)$ is larger than the space
of square integrable functions with respect to any smooth Riemannian metric on $S$, by the
Poincar\'e inequality the space $H^1_\omega(S)$ coincides as a topological vector space with
the standard Sobolev space $H^1(S)$ on the closed manifold $S$. In particular,  it does not depend on
the abelian differential $\omega\in \Ab$. However the natural norm on $H^1_\omega(S)$, defined for all $f\in H^1_\omega(S)$ as
\begin{equation}
\vert f  \vert_{1,\omega}:=  \left( \vert f \vert_{0, \omega}^2 + \vert X_\omega f\vert_{0,\omega}^2 + 
\vert Y_\omega f \vert_{0,\omega}^2 \right)^{1/2} \,,
\end{equation}
does depend on $\omega\in \Ab$ and it is not invariant under the $SL(2, \R)$ or
even the $g_t$-action. In fact, the following optimal bound holds. For any $f \in H^1(S)$ (we recall
that $H^1(S)$ and $H^1_\omega(S)$ coincide as topological vector spaces),
\begin{equation}
\label{eq:normdist}
\vert f \vert_{1, g_t\omega} \leq   e^{\vert t \vert} \, \vert f \vert_{1,\omega} \,,  \quad \text{ \rm for any }\,
t\in \R  \,.
\end{equation}

 The weighted Sobolev spaces $W^1_\omega(S)$ of $1$-forms on $S$ are defined as follows. Let $\omega \in \Ab$ and let $\Psi^1(S\setminus \Sigma_\omega)$ be the space of all measurable $1$-forms on $S\setminus \Sigma_\omega$. We define
\begin{equation}
W^1_\omega(S) := \{\psi\in \Psi^1(S\setminus \Sigma_\omega)  \vert  (\imath_{X_\omega} \psi , \imath_{Y_\omega} \psi) \in H^1_\omega(S)
\times  H^1_\omega(S)\}\,.
\end{equation}
By definition, the Hilbert space $W^1_\omega(S) \equiv H^1_\omega(S) \times  H^1_\omega(S)$, 
hence the topological vector space underlying the Hilbert space $W^1_\omega(S)$ is independent of $\omega\in \Ab$ up to isomorphisms and, for any stratum $\Omega \subset \Omega_g$, it is independent of $\omega \in \Omega$. In particular, as a topological vector space $W^1_\omega(S)$ is invariant under the action of $SL(2, \R)$ on $\omega \in \Ab$ (while the Hilbert structure is not). 

 For any $\omega\in \Ab$, we will adopt the standard notation for the 
dual Hilbert spaces, that is,
\begin{equation}
H^{-1}_\omega(S):= H^1_\omega(S)^\ast \quad \text{ and } \quad  
W^{-1}_\omega(S) :=W^1_\omega(S)^\ast \,.
\end{equation}
Let $\vert \cdot \vert_{-1, \omega}$ denote the Sobolev norm on $H^{-1}_\omega(S)$ or $W^{-1}_\omega(S)$. We remark that the topological vector space underlying the dual space 
$H^{-1}_\omega(S)$ coincides with the standard dual Sobolev space $H^{-1}(S)$. 

 \subsection{Distributional cocycles}
\label{distcocycles}
We recall the definition, first introduced in \cite{Forni}, of a natural lift of the Kontsevich-Zorich 
cocycle to a bundle of currents over the moduli space $\Omega_g$ of abelian differentials.

 The identity component $\text{\rm Diff}^+_0(S)$ of the group $\text{\rm Diff}^+(S)$ of orientation-preserving diffeomorphisms of the oriented surface $S$ acts in a natural way on the trivial bundle $\Ab \times H^{-1}(S)$. In fact, any $\phi\in \text{\rm Diff}^+(S)$ defines by pull-back an isomorphism $\phi^* :  H^{-1}(S) \to   H^{-1}(S)$. 
 The quotient bundle
\begin{equation}
\label{eq:distbundleone}
 \HTeich^{-1}_g(S) := \Ab \times H^{-1}(S) / 
\text{\rm Diff}^+_0(S)
\end{equation}
is a well-defined orbifold vector bundle over the Teichm\"uller space $\Teich$ of abelian holomorphic differentials. There is natural action of the mapping class group $\Gamma_g$ on the bundle \pref{eq:distbundleone} induced by the action of $\text{\rm Diff}^+(S)$ on the trivial bundle $ \Ab \times H^{-1}(S)$. The resulting quotient bundle
\begin{equation}
\label{eq:distbundletwo}
H^{-1}_g(S):= \HTeich^{-1}_g(S)/\Gamma_g \,\,,\\
\end{equation}
is a well-defined orbifold vector bundle over the moduli space $\Omega_g$ of
abelian differentials. Since the Hilbert norms on the spaces
$H^{-1}_\omega(S)$ are equivariant under the action of the group $\text{\rm Diff}^+(S)$ 
on $\omega\in \Ab$, the bundle $\HTeich^{-1}_g(S)$ and its
quotient $H^{-1}_g(S)$ have natural structures of (orbifold) Hilbert bundles.

 Let $\{G_t \vert t\in\R\}$ be the cocycle over the Teichm\"uller flow  on $\Omega_g$ defined  as the projection onto the bundle $H^{-1}_g(S)$ of the trivial cocycle 
\begin{equation}
\label{eq:distcocycle1}
 g_t \times\text{\rm id} :  \Ab \times H^{-1}(S)  \to \Ab \times H^{-1}(S)\,\,.
\end{equation}

 A similar construction allows us to define a cocycle over the Teichm\"uller flow on a Sobolev bundle of $1$-dimensional currents. It can be proved that the space
$$
{\WAb}^{-1}_g(S):= \bigcup_{\omega\in \Ab}  
\{\omega\} \times W^{-1}_\omega(S) 
$$
has the structure of a Hilbert bundle over $\Ab$ on which the group
$\text{\rm Diff}^+(S)$ acts naturally by pull-back. The quotient bundles
\begin{equation}
\label{eq:distbundle2}
\begin{aligned}
 \WTeich^{-1}_g(S) &:=    {\WAb}^{-1}_g(S) / \text{\rm Diff}^+_0(S)\,,\\ 
 W^{-1}_g(S) &:=   \WTeich^{-1}_g(S) / \Gamma_g\,\,
\end{aligned}
\end{equation}
are well-defined orbifold vector bundles over the Teichm\"uller space $\Teich$ and
over the moduli space $\Omega_g$ respectively. Since the Hilbert norms on the spaces
$W^{-1}_\omega(S)$ are equivariant under the action of $\text{\rm Diff}^+(S)$ 
on $\omega\in \Ab$, the bundle $\WTeich^{-1}_g(S)$ and its
quotient $W^{-1}_g(S)$ have natural structures of (orbifold) Hilbert bundles.

  Let  $\{{\dc}_t\vert t\in \R\}$ be the cocycle over the Teichm\"uller flow on $\Omega_g$ 
defined as the projection onto the bundle $W^{-1}_g(S)$  of the trivial (skew-product) cocycle
given for each $(\omega,t) \in \Ab\times \R$ by the identity map
\begin{equation}
\label{eq:distcocycle2}
  \text{\rm id} :  W^{-1}_\omega(S)  \to W^{-1}_{g_t\omega}(S)\,\,.
\end{equation}
The identity maps in \pref{eq:distcocycle2} are well-defined since, as remarked above, the topological vector space $W^{-1}_{\omega}(S)$ is invariant under the action of $SL(2,\R)$, hence of the Teichm\"uller flow $\{g_t \vert t\in\R\}$,  on $\omega\in \Ab$.

 The above cocycles can be defined in terms of the parallel transport 
of distributions and $1$-currents with respect to the standard (Gauss-Manin) flat
connections on the bundles $H^{-1}_g(S)$ and $W^{-1}_g(S)$ along the orbits of the 
Teichm\"uller flow.

By the definition of the Teichm\"uller flow and of the distributional cocycles, the following isomorphism holds:
\begin{equation}
\label{eq:currcocycle}
{\dc}_t \equiv \text{\rm diag}(e^{-t},e^t)\otimes G_t \quad \text{ \rm on } \quad 
W_g^{-1}(S)\equiv \R^2\otimes H^{-1}_g(S)\,.
\end{equation}
 By \pref{eq:normdist} and \pref{eq:currcocycle} we can immediately derive
that, for any $(\omega,t)\in \Omega_g\times \R$, 
\begin{equation}
\label{eq:cocycleestimates}
\begin{aligned}
\vert G_t (\Cal D) \vert_{-1, g_t\omega}   &\leq    e^{\vert t \vert }\, \vert \Cal D \vert_{-1,\omega}
\, , \quad \text{ for any } \, \Cal D \in H^{-1}_\omega(S) \,; \\
\vert {\dc}_t (\gamma) \vert_{-1, g_t\omega}   &\leq    e^{2 \vert t \vert }\, \vert \gamma \vert_{-1,\omega}
\, , \quad \text{ for any } \, \gamma \in W^{-1}_\omega(S) \,.
\end{aligned}
\end{equation} 
In the following we will essentially be concerned with the restriction of the cocycle $\{{\dc}_t\vert t\in\R\}$ to the invariant  Hilbert sub-bundle ${\Cal Z}^{-1}_g(S)\subset W^{-1}_g(S)$ of closed currents. By definition, the fiber of the bundle ${\Cal Z}^{-1}_g(S)$ is given at any $\omega\in \Omega_g$ by the subspace of \emph{closed currents}:
\begin{equation}
\label{eq:Zomega}
{\Cal Z}^{-1}_\omega(S) := \{ \gamma\in W^{-1}_\omega(S)\vert d\gamma=0 \}\,.
\end{equation}
The exterior derivative operator $d$ is defined in the weak sense with respect to the appropriate space of test functions (see \S 6 in \cite{Forni}). By Lemma 6.2 in \cite{Forni}, the generalized de Rham theorem implies that for each $\omega\in \Ab$ there is a linear surjective cohomology map 
\begin{equation}
\label{eq:jomega}
j_\omega : {\Cal Z}^{-1}_\omega(S) \to H^1(S,\R)\,,
\end{equation}
hence there is a bundle \emph{cohomology map }
\begin{equation}
\label{eq:j}
j : {\Cal Z}^{-1}_g(S) \to {\Cal H}^1_g(S,\R)
\end{equation}
onto the cohomology bundle ${\Cal H}^1_g(S,\R)$ over the moduli space $\Omega_g$. The kernel ${\Cal E}_g^{-1}(S) \subset {\Cal Z}^{-1}_g(S)$ of the cohomology map consists of the smooth Hilbert sub-bundle of \emph{exact currents}, which can be described as follows. For any $\omega\in \Omega_g$, let
\begin{equation}
\label{eq:Eomega}
{\Cal E}^{-1}_\omega(S) := \{ \gamma \in W^{-1}_\omega(S)\vert  \gamma=dU \,,  \quad U \in L^2_\omega(S)  \}\,.
\end{equation}
It follows from \cite{Forni}, Lemma 7.6, that the space ${\Cal E}^{-1}_\omega(S)$ does coincide with
the kernel of the cohomology map \pref{eq:jomega} for all $\omega\in \Omega_g$. 

For any stratum $\Omega \subset \Omega_g$, let ${\Cal Z}^{-1}_\Omega(S)$ and 
${\Cal E}^{-1}_\Omega(S)$ be the restrictions to $\Omega$ of the bundles ${\Cal Z}^{-1}_g(S)$ and 
${\Cal E}^{-1}_g(S)$ respectively and let ${\Cal H}_\Omega^1(S,\R)$ be the restriction to $\Omega$
of the cohomology bundle ${\Cal H}_g^1(S,\R)$. The above discussion leads to the following:

\begin{lemma}  \label{lm:deRhamth} For any stratum $\Omega\subset \Omega_g$ of the moduli space, the bundle cohomology map \pref{eq:j} induces a smooth isomorphism of smooth normed orbifold vector bundles
\begin{equation}
i_\Omega:  {\Cal Z}^{-1}_\Omega(S)/{\Cal E}^{-1}_\Omega(S) \,\, \to \,\, {\Cal H}_\Omega^1(S,\R)\,.
\end{equation}
\end{lemma}
\begin{proof} 
Let $\hat \Omega \subset \hat \Omega_g$ be the preimage of $\Omega$ with respect to the natural projection $\hat \Omega_g \to \Omega_g= \hat \Omega_g/\text{\rm Diff}^+(S)$, consisting of all abelian differentials with zeroes at a fixed set $\Sigma_\Omega$ of the prescribed multiplicities. The set $\hat\Omega$ endowed with the smooth topology has the structure of a Fr\'echet space. 

The pull-back of the cohomology bundle ${\Cal H}_\Omega^1(S,\R)$ to $\hat \Omega$ is by definition locally isomorphic to the trivial bundle $\hat \Omega \times H^1(S, \R)$. The pull-back of the de Rham bundle ${\Cal Z}^{-1}_\Omega(S)/{\Cal E}^{-1}_\Omega(S)$ to $\hat \Omega$ can be smoothly identified to $\hat \Omega \times H^1(S, \R)$. In fact, by definition of the Sobolev norms on the spaces of $1$-forms and $1$-currents, the topological vector spaces ${\Cal E}^{-1}_\omega(S) \subset {\Cal Z}^{-1}_\omega(S)$ are  independent of the basepoint $\omega\in \hat \Omega$ up to an isomorphism which depends smoothly on the basepoint. With respect to the local smooth trivializations just described, the map $i_\Omega$ is the projection of the identity map, hence it is a smooth isomorphism of orbifold vector bundles.

For any $\omega\in \hat \Omega_g$ the finite dimensional vector space $H^1(S,\R)$ can be endowed with the following (equivalent) norms: the Sobolev norm $\vert \cdot \vert_{-1,\omega}$, defined as the quotient norm on the quotient space $\Cal Z^{-1}_\omega(S)/\Cal E^{-1}_\omega(S)$ and the Hodge norm $\Vert \cdot \Vert_\omega$, introduced in \pref{eq:Hodgenorm}. It can be proved that both these norms depend smoothly on the differential $\omega\in \hat \Omega$, hence ${\Cal Z}^{-1}_\Omega(S)/{\Cal E}^{-1}_\Omega(S)$ and ${\Cal H}_\Omega^1(S,\R)$ are smooth normed orbifold vector bundles.
\end{proof}

It follows from Lemma \ref{lm:deRhamth} that there exists a strictly positive continuous function $K_\Omega: \Omega \to \R^+$ such that, for all abelian differentials $\omega\in \Omega$ and all cohomology classes $c\in H^1(S,\R)$,
\begin{equation}
\label{eq:normscomp}
K_\Omega (\omega) \, \Vert c \Vert _\omega \leq    
\vert c \vert _{-1,\omega} \leq \Vert c \Vert _\omega\,.
\end{equation}
In fact, the first inequality in \pref{eq:normscomp} holds since the bundle cohomology map induces an isomorphism of topological normed bundles. The second inequality in \pref{eq:normscomp} holds by the definition of the quotient norm as an infimum over each equivalence class and by the definition \pref{eq:Hodgenorm} of the Hodge norm. In fact, for any $\omega \in \hat \Omega$, the Hodge norm of a cohomology class $c\in H^1(S,\R)$ is the norm $\vert m \vert_{0,\omega}$ of a meromorphic function $m\in L^2_\omega(S) \subset H^{-1}_\omega(S)$ such that $\operatorname{Re}(m\,\omega)$ is a harmonic representative. By its definition and a simple computation, the norm of the $1$-form $\operatorname{Re}(m\,\omega)$ as a current in $W^{-1}_\omega(S)$ is equal to $\vert m \vert_{-1,\omega} \leq \vert m \vert_{0,\omega}$, hence the second inequality in \pref{eq:normscomp} holds.

The analysis of the Lyapunov structure of the cocycle ${\dc}_t \vert {\Cal Z}^{-1}_g(S)$
(carried out in \cite{Forni} with respect to $g_t$-invariant probability measures which are non-uniformly hyperbolic for the Kontsevich-Zorich cocycle) is based on the following facts. The cocycle projects
onto the Kontsevich-Zorich cocycle on the cohomology bundle $ {\Cal H}^1_g(S,\R)$, in the sense that
\begin{equation}
j \circ {\dc}_t   =  \Phi_t \circ j    \quad \text{ on } \,\, {\Cal Z}^{-1}_g(S) \,.
\end{equation}
In addition, the restriction ${\dc}_t \vert {\Cal E}^{-1}_g(S)$ to the sub-bundle of exact currents
has a single (well-defined) Lyapunov exponent equal to $0$. In fact, a stronger statement holds:
there exists a ${\dc}_t$-invariant norm on  ${\Cal E}^{-1}_g(S)$ which is equivalent to the
Sobolev norm $\vert \cdot \vert_{-1}$ on ${\Cal E}^{-1}_g(S)\subset W^{-1}_g(S)$ on compact sets.

 For any $\gamma\in {\Cal E}^{-1}_\omega(S)$, there exists a unique function $U_\gamma \in L^2_\omega(S)$ such that 
\begin{equation}
\gamma = dU_\gamma    \quad \text{and}  \quad   \int_S U_\gamma\,dA_\omega \,\,=\,\,0 \,\,.
\end{equation}

\begin{lemma} (Lemma 9.3 in \cite{Forni}) \label{lm:Lnormineq} 
The norm $\Vert \cdot \Vert_{-1}$, defined on each fiber  ${\Cal E}^{-1}_\omega(S)$ of the bundle 
${\Cal E}^{-1}_g(S)$ of exact currents as
\begin{equation}
\label{eq:Lyapunovnorm}
\Vert \gamma\Vert_{-1,\omega} :=\vert U_ \gamma \vert_{0,\omega}\,\,,\quad \text{ for any } \, 
 \gamma \in {\Cal E}^{-1}_\omega(S)\,
\end{equation}
is invariant under the distributional cocycle $\{{\dc}_t \vert t\in \R\}$. There exists a continuous function 
$K_g:\Omega_g \to \R^+$ such that, for all $\omega \in \Omega_g$,
\begin{equation}
\label{eq:Lnormineq} 
K_g(\omega)\, \Vert  \gamma \Vert_{-1,\omega} \leq \vert  \gamma \vert_{-1,\omega} \leq 
\Vert  \gamma \Vert_{-1,\omega}\,\,,\quad \text{ for all }  \gamma \in {\Cal E}^{-1}_\omega(S)  \,\,.
\end{equation}
\end{lemma}
 We remark that the invariance of the norm \pref{eq:Lyapunovnorm} claimed above follows immediately from the invariance of the norm $\vert \cdot \vert_{0,\omega}$ on $L^2_\omega(S)$ under the action of the Teichm\"uller flow $g_t$ on $\omega\in \Omega_g$.

\section{Spectral Gap}
\label{spectral}
In this section we prove a spectral gap theorem for the Kontsevich-Zorich coycle on the cohomology
bundle which is then extended to the restriction of the distributional cocycle to appropriate closed invariant subsets of the bundle of currents, which contain all currents given by integration along curves. The gap theorem for the distributional cocycle is a key ingredient in the proof of our main theorem
on the decay of ergodic averages. 

\subsection{The Kontsevich-Zorich cocycle}
\label{KZgap}
 
 Let $\omega\in \Omega_g$ be an abelian differential and let $c\in H^1(M_\omega,{\R})$ be a cohomology class. Let $\omega_t = g_t \omega$, let $M_t=M_{\omega_t}$ be the underlying Riemann surface and let $\Cal M_t = \Cal M_{\omega_t} \subset L^2_\omega(S)$ be the subspace of meromorphic functions on $M_t$.  By \pref{eq:repr} there exists a one-parameter family $\{m_t\}_{t\in \R} \subset  \Cal M_t$ (implicitly) defined by the identity
\begin{equation}
\label{eq:ct}
\Phi_t(c) =c_{\omega_t} (m_t):=[\operatorname{Re}(m_t \,\omega_t)] \in H^1(M_t,{\R})\,.
\end{equation}

Let $B_{\omega}:L^{2}_{\omega}(S)\times  L^{2}_{\omega}(S) \to \C$ be the complex bilinear form given by
\begin{equation}
\label{eq:Bform}
B_{\omega}(u,v) := \int_{M} u \,v \, dA_{\omega}\, ,\quad \hbox{ for all } u,v \in L^{2}_{\omega}(M)\,.
\end{equation}
\begin{lemma}  
\label{lemma:normder} (see \cite{Forni}, Lemma 2.1') The variation of the Hodge
norm $\Vert \Phi_t(c) \Vert_{\omega_t}$, which coincides with the $L^{2}_{\omega}$-norm $\vert m_{t}\vert_{0}$
under the identification \pref{eq:ct}, is given by the following formula:
\begin{equation}
\label{eq:normder} 
\frac{d}{dt} \vert m_t \vert_0^2=-2\,\operatorname{Re}\,B_{\omega}(m_t)  = -2\, \operatorname{Re} \,\int_S (m_t)^2 dA_{\omega}
\end{equation}
\end{lemma}
The following bound on the Hodge norm of the Kontsevich-Zorich cocycle can be easily derived from Lemma \ref{lemma:normder}. For any $\omega \in \Omega_g$ and all $t\in \R$, 
\begin{equation}
\label{eq:KZnorm1}
\Vert \Phi_t: H^1(M_\omega, \R) \to  H^1(M_{\omega_t}, \R) \Vert \,\, \leq \,\,   e^{\vert t \vert }  \,\,.
\end{equation}
In fact, it follows from formula \pref{eq:normder} by the Cauchy-Schwarz inequality that for any
$c\in H^1(M_\omega, \R)$ and all $t\in \R$, 
\begin{equation}
\frac{d}{dt} \Vert \Phi_t (c)  \Vert_{\omega_t} ^2 \leq 2\, \Vert \Phi_t (c)  \Vert_{\omega_t} ^2\,.
\end{equation}

 We prove below a spectral gap theorem for Lebesgue almost all differentials in {\it any }orbit of the circle group $SO(2,\R)$ on any stratum $\Omega\subset \Omega_g$. The argument is based on formula \pref{eq:normder} and on Theorem~\ref{proportion}.

\smallskip
Let $I(M_\omega,\R) \subset H^1(M_\omega, \R)$ be the $2$-dimensional subspace defined as 
\begin{equation}
\label{eq:T}
I(M_\omega,\R):= \R \cdot \operatorname{Re}(\omega) + \R \cdot \operatorname{Im}(\omega) 
\end{equation}
and let $I^\perp (M_\omega,\R)$ be the symplectic orthogonal of $I(M_{\omega},\R)$ in $H^1(M_\omega, \R)$, with respect to the symplectic structure induced by the intersection form:
\begin{equation}
\label{eq:Tperp}
I^\perp (M_\omega,\R):=\{ c\in H^1(M_\omega, \R)   \vert    c \wedge [\omega] =0\}\,.
\end{equation}
 The complementary sub-bundles $I_g(S,\R)$ and $I^\perp_g(S,\R) \subset \Cal H^1_g(S,\R)$, with fibers at any $\omega\in \Omega$ respectively equal to $I(M_\omega,\R)$ and $I^\perp(M_\omega,\R)$, are invariant under the Kontsevich-Zorich cocycle. In fact, it is immediate to verify that the sub-bundle $I_g(S,\R)$ is invariant under the Kontsevich-Zorich cocycle and that the Lyapunov spectrum of the restriction of the Kontsevich-Zorich cocycle to $I_g(S,\R)$ equals $\{1,-1\}$ (both exponents with multiplicity $1$). Hence, taking into account the upper bound \pref{eq:KZnorm1}, it follows that the top Lyapunov exponent of the Kontsevich-Zorich cocycle is equal to $1$. The invariance of the symplectic orthogonal bundle $I^\perp_g(S,\R)$ follows since the Kontsevich-Zorich cocycle is symplectic.

\begin{lemma} 
\label{lm:KZnorm2}  There exists a continuous function $\Lambda: \Omega \to [0, 1)$ such that
for all $\omega \in \Omega$ and all $t\in \R$, the following bound holds:
\begin{equation}
\label{eq:KZnorm2}
\Vert \Phi_t: I^\perp(M_\omega, \R) \to  I^\perp(M_{\omega_t}, \R) \Vert \,\, \leq \,\,   \exp \left(
\int_0^{\vert t \vert} \Lambda (\omega_s)\, ds  \right)\,. 
\end{equation}
\end{lemma}

\begin{proof} The argument follows closely the proof of Corollary 2.2 in \cite{Forni}. 

By formula
\pref{eq:normder},
\begin{equation}
\label{eq:derlog1}
\frac{d}{dt} \log \vert m_t \vert_0=- \, { \frac{\operatorname{Re} [B_{\omega}(m_t)]}{\vert m_t\vert_0^2}}\,.
\end{equation}
Under the isomorphism \pref{eq:repr}, the subspace $I^\perp(M_\omega,\R)$ is represented  by meromorphic functions with {\it zero average }(orthogonal to constant functions). Hence, following \cite{Forni}, we define a  function $\Lambda:  \Omega \to \R^+$ as follows: for any $\omega\in  \Omega$, 
\begin{equation}
\label{eq:Lambdaplus}
\Lambda(\omega):= \max \{ { \frac{\vert B_{\omega}(m)\vert}{\vert m\vert_0^2}}\,|\, m\in {\Cal M}_{\omega}\setminus\{0\}\,,\,\,\int_S m\,dA_{\omega}=0\,\} \,. 
\end{equation}
By  \pref{eq:derlog1} and \pref{eq:Lambdaplus} it follows that for any $c \in I^\perp(M_\omega, \R)$ and all $t\in \R$,
\begin{equation}
\label{eq:derlog2}
\frac{d}{dt} \log \Vert \Phi_t(c) \Vert_{\omega_t} \leq    \Lambda(\omega_t) \,.
\end{equation}
and the desired upper bound \pref{eq:KZnorm2} follows by integrating formula \pref{eq:derlog2}.

The function $\Lambda:  \Omega \to \R^+$  is continous since the  Hilbert space $L^2_\omega(S)$ 
and the finite dimensional  subspace $\Cal M_\omega \subset L^2_\omega(S)$ of meromorphic functions depend continuously on $\omega \in \Omega$.  

It remains to be proven that the function $\Lambda(\omega) <1$ for all $\omega\in \Omega$.
Since by the Cauchy-Schwarz inequality, for any $m\in\Cal M_\omega$, 
\begin{equation}
\label{eq:Schwarz}
\vert B_{\omega}(m)\vert=\vert(m, {\overline{m}})_{\omega} \vert \leq \vert m \vert_0^2\,\,, 
\end{equation}
the range of the function $\Lambda$ is contained in the interval $[0,1]$. We claim
that $\Lambda(\omega)<1$ for all $\omega\in \Omega$.  In fact, $\Lambda(\omega)=1$ if 
and only if there exists a {\it non-zero }meromorphic function with zero average $m\in {\Cal M}_{\omega}$ such that $\vert(m,{\overline {m}})_{\omega}\vert=\vert m\vert_0^2$. A well-known property of the Cauchy-Schwarz inequality then implies that there exists $u \in {\C}$ such that $m=u \,{\overline{m}}$. However, it cannot be so, since $m$ would be meromorphic and anti-meromorphic, hence constant, and by the zero average condition it would be zero.
\end{proof}

\begin{lemma}
\label{lm:recurrest}
Let $\Lambda : \Omega \to [0,1)$ be a continuous function. There exists $\lambda \geq 0$ such that,
for any $\omega \in \Omega$ and for almost all $\theta\in \Sone$,
\begin{equation}
\label{eq:recurrest}
\limsup _{t\to +\infty} \frac{1}{t}  \int_0^t   \Lambda (g_s r_\theta \omega) \, ds \,\, \leq \,\, \lambda
<1\,.
\end{equation}
\end{lemma}
\begin{proof} By the large deviations result Theorem~\ref{proportion},  for any strictly positive number $0<\eta<1$ there exists a compact set $C=C(\eta) \subset  \Omega$ such that, for all $\omega \in \Omega$ and for almost all $\theta\in \Sone$, the following holds:
\begin{equation}
\label{eq:largedev}
\limsup_{t\to \infty} \frac{1}{t} \, \vert\{ 0 \leq s \leq t\, \vert \, g_s r_\theta \omega \not\in C \} \vert
\,\,  \leq \,\, \eta \,.
\end{equation}
Let $\Lambda_C:= \max \{ \Lambda(\omega) \,\vert \, \omega\in C\}$ and, for any $(t,\omega)\in \R^+\times \Omega$, let 
$$
\Cal E_C(t,\omega):= \text{ \rm Leb}\left (\{ 0 \leq s \leq t\, \vert \, g_s \omega \not\in C \} \right)\,.
$$ 
Since the function $\Lambda$ is continuous and $\Lambda(\omega)<1$ for all $\omega\in \Omega$, its maximum on any compact set is strictly less than $1$, in particular $\Lambda_C <1$. The following immediate inequality holds:
\begin{equation}
\label{eq:logintsplit}
\int_0^t  \Lambda ( g_s r_\theta \omega )\, ds  \leq  
(1 -  \Lambda_C) \,  {\Cal E}_C(t,r_\theta\omega) + t \, \Lambda_C\,.
\end{equation}
It follows from the above argument that the number 
$$
\lambda=(1 -  \Lambda_C) \eta \,+\, \Lambda_C <1
$$ 
can be chosen independently of the abelian differential $\omega \in \Omega$.
\end{proof}

We define the {\it upper second exponent }of the Kontsevich-Zorich cocycle at any differential $\omega\in\Omega$ as the top upper Lyapunov exponent at $\omega$ of the restriction of the cocycle to the sub-bundle $I^\perp_g(S,\R)$, that is the number 

\begin{equation}
\label{eq:l2plus}
\lambda^+_2(\omega) :=  \limsup_{t\to +\infty}  \frac{1}{t} \log  \Vert
 \Phi_t \vert I^\perp(M_{\omega},\R) \Vert \,.
\end{equation}

\begin{theorem}
\label{thm:spectralgap}
 \text{\rm (Spectral gap)} For any stratum $\Omega$, there is a (non-negative) real number 
 $\lambda=\lambda(\Omega)$ such that for any  $\omega\in \Omega$, 
\begin{equation}
\label{eq:specgap}
 \lambda^+_2(r_\theta \omega)   \, \leq  \, \lambda\,<\, 1\,,  \quad \text {\rm for almost all } 
 \, \,\theta\in \Sone\,.
\end{equation}
\end{theorem}

\begin{proof}
By Lemma \ref{lm:KZnorm2} there exists a continuous function $\Lambda: \Omega
 \to [0,1)$ such that for any $\omega\in \Omega$ and any $t\in \R^+$,
\begin{equation}
\label{eq:logbound}
 \frac{1}{t} \log  \Vert  \Phi_t \vert I^\perp(M_\omega,\R) \Vert \leq   \frac{1}{t}  \int_0^t  
 \Lambda( g_s\omega )\, ds \,.
\end{equation}
The statement then follows immediately from Lemma \ref{lm:recurrest}.
\end{proof}

\subsection{The distributional cocycle}
\label{ss:distrgap}

We prove below a spectral gap theorem for the cocycle $\{{\dc}_t \vert t\geq 0\}$ on the bundle of currents.

For each $\omega\in \Omega$, we let 
\begin{equation}
\begin{aligned}
\Cal I_\omega(S) &=  \R \operatorname{Re} (\omega)  \,\oplus\,  \R   \operatorname{Im} (\omega) \subset  W^{-1}_\omega(S) \,; \\
\Cal I^\perp_\omega(S)&= \{ \gamma \in W^{-1}_\omega(S)\, \vert \, \<\gamma \wedge \omega,1 \>=0\}\,.
\end{aligned}
\end{equation}
Let $\Cal I_g(S)$ and  $\Cal I^\perp_g(S) \subset W^{-1}_g(S)$ be the sub-bundles with fibers
equal to $\Cal I_\omega(S)$ and $\Cal I^\perp_\omega(S)$ at any $\omega\in \Omega$.  There
is a ${\dc}_t$-invariant splitting
\begin{equation}
\label{eq:Wsplit}
W^{-1}_g(S)= \Cal I_g(S)    \oplus  \Cal I^\perp_g(S) \,.
\end{equation}
Let $\Cal Z^{-1}_g(S) \subset W^{-1}_g(S)$ be the sub-bundle of closed currents.  The bundles  $\Cal I_g(S) \subset W^{-1}_g(S)$ and $\Cal I^\perp_g(S) \cap W^{-1}_g(S)$ project onto the complementary sub-bundles $I_g(S,\R)$ and  $I^\perp_g(S,\R) \subset \Cal H_g(S,\R)$, respectively, under the cohomology map. The splitting \pref{eq:Wsplit} induces a ${\dc}_t$-invariant splitting
\begin{equation}
\label{eq:Zsplit}
\Cal Z^{-1}_g(S)= \Cal I_g(S)    \oplus \left( \Cal I^\perp_g(S)\cap  \Cal Z^{-1}_g(S) \right)\,.
\end{equation}
Let $\delta_g : W^{-1}_g(S) \to \R$ be the (continuous) distance functions to the Hilbert sub-bundle
$\Cal Z^{-1}_g(S)$ of closed currents defined as follows:  for each $\omega \in \Omega$, the restriction $\delta_g \vert  W^{-1}_\omega(S)$ is equal to the distance function from
the closed subspace $\Cal Z^{-1}_\omega(S) \subset W^{-1}_\omega(S)$ with respect to the Hilbert
space metric on $W^{-1}_\omega(S)$.  We introduce the following closed, ${\dc}_t$-invariant subsets $\Gamma_C(\delta)$ of the bundle $W^{-1}_g(S)$. For any compact set $C \subset \Omega$ 
and any $\delta >0$, let $\Gamma_C(\delta) \subset W^{-1}_g(S)$ be the set defined as follows:
\begin{equation}
\label{eq:GammaCdelta}
\Gamma_C(\delta) \cap W^{-1}_\omega(S) =\{ \gamma \in  W^{-1}_\omega(S) \,\vert \,
g_t \omega \in C \Rightarrow  \delta_g \left( {\dc}_t (\gamma) \right) \leq \delta \} \,.
\end{equation}
In other terms, the fibered subset $\Gamma_C(\delta)$ contains all currents which stay at bounded distance ($\leq \delta$) from the sub-bundle of closed currents for all returns of the Teichm\"uller orbit
to a given compact set $C\subset \Omega$. The relevant examples of non-closed currents in
$\Gamma_C(\delta)$ are given by currents of integration along paths in $S$. In fact, as we shall see below, for any compact set $C\subset \Omega$ there exists $\delta_C >0$ such that any  current represented by a path on the surface $S$ belongs to $\Gamma_C(\delta)$ for $\delta\geq \delta_C$.

The core technical result of this paper is the following `spectral gap' theorem for the restriction of
the distributional cocycle $\{{\dc}_t \vert t\in \R\}$ to any invariant set $\Gamma_C(\delta)
\subset W^{-1}_g(S)$.

We introduce the following notation. For any $\omega\in\Omega$, let $t_0=0$ and let $\{t_n \vert n\in\Z^+\}$ denote a non-decreasing sequence of visiting times of the forward orbit $\{g_t\omega \vert t \geq 0\}$ to a given compact set $C\subset \Omega$. For all $n\in \Z^+$, let $\omega_n= g_{t_n} \omega\in C$ and, for any $\gamma\in W^{-1}_\omega(S)$, let $\gamma_n= \dc_{t_n}(\gamma) \in W^{-1}_{\omega_n}(S)$. 

For any $r\in \R$, let us denote
$$
r^+ =   \max \{ 1, r \}\,.
$$ 
Let $\Lambda:\Omega \to [0,1)$ be the continous function defined in Lemma \ref{lm:KZnorm2}.

\begin{lemma} 
\label{lm:iterativest}
For any compact set $C\subset \Omega$, there exists a constant $K_C>1$ such that,
for any $\omega\in \Omega$, for any $\gamma \in \Cal I^\perp_\omega(S)\cap \Gamma_C(\delta)$ and for all $k\in \Z^+$, the following estimate holds:
\begin{equation}
\label{eq:iterativest}
\vert \gamma \vert_{-1, \omega} \leq K_C \, \delta^+
 \vert  \gamma_k \vert^+_{-1, \omega_k}  \exp\bigl( \int_0^{t_k} \Lambda(g_s\omega) \,ds\bigr)
   [ \sum_{j=0}^{k-1}  e^{2(t_{j+1}-t_j)} ]^3\,.
 \end{equation}
\end{lemma}

\begin{proof} For each $j \in \Z^+$, since $ \Cal Z^{-1}_{\omega_j} (S)$ is closed in 
$W^{-1}_{\omega_j} (S)$, there exists an orthogonal decomposition,
\begin{equation}
\label{eq:gammasplit}
 {\dc}_{t_j} (\gamma) = z_j \, + \, r_j  \,, \quad  \text{ \rm with } \, z_j \in \Cal Z^{-1}_{\omega_j}(S) \,,\, \, r_j \perp \Cal Z^{-1}_{\omega_j}(S)\,,
\end{equation}
and, since $\gamma  \in \Gamma_C(\delta)$ and $\omega_j\in C$, the following crucial bound holds:
\begin{equation}
\label{eq:Gammadeltabound}
\vert r_j \vert_{-1, \omega_j} \leq  \delta \,.
\end{equation}
For each $j\in \Z^+$, let $\pi_j: W^{-1}_{\omega_j} (S) \to \Cal Z^{-1}_{\omega_j} (S)$ denote the orthogonal projection  and let $\tau_j = t_{j+1} -t_{j}$. By \pref{eq:gammasplit} and by orthogonal projection on the ${\dc}_t$-invariant bundle $\Cal Z^{-1}_g(S)$ the following recursive identity holds:
\begin{equation}
\label{eq:recid1}
z_{j} = {\dc}_{-\tau_j} (z_{j+1}) \, + \,    \pi_{j}\circ  {\dc}_{-\tau_j}(r_{j+1})  \in 
\Cal Z^{-1}_{\omega_j}(S)\,.
\end{equation}
By the bounds  \pref{eq:cocycleestimates} and  \pref{eq:Gammadeltabound}, it follows that
\begin{equation}
\label{eq:remainderbound1}
\vert \pi_{j}\circ  {\dc}_{-\tau_j}(r_{j+1}) \vert_{-1, \omega_j}  \leq 
\vert {\dc}_{-\tau_j}(r_{j+1}) \vert_{-1, \omega_j}  \leq e^{2\tau_j} \, \delta \,.
\end{equation}

By projection on the cohomology bundle $\Cal H^1_g(S,\R)$ and by the comparison estimate
\pref{eq:normscomp}, we derive from the identity \pref{eq:recid1} and from the bound \pref{eq:remainderbound1} that there exists $K^{(1)}_C >1$ such that 
the Hodge norm
\begin{equation}
\label{eq:recbound1}
\Vert [z_j] - \Phi_{-\tau_j} ([z_{j+1}]) \Vert_{\omega_j}  \leq K^{(1)}_C \, \delta\,e^{2\tau_j} \,.
\end{equation}
For all $j\in\Z^+$, since $\gamma \in {\Cal I}^\perp_\omega(S)$ and the sub-bundle ${\Cal I}^\perp_g(S)$ is ${\dc}_t$-invariant, ${\dc}_{t_j} (\gamma) \in {\Cal I}^\perp_{\omega_j}(S)$. We also
have that $r_j \in {\Cal I}^\perp_{\omega_j}(S)$ since $\Cal I_g(S) \subset \Cal Z^{-1}_g(S)$ and $r_j\in W^{-1}_{\omega_j}(S)$ is orthogonal to $\Cal Z^{-1}_{\omega_j}(S)$ by definition. Hence $z_j \in {\Cal I}^\perp_{\omega_j}(S)$ and its de Rham cohomology class $[z_j] \in I^\perp_{\omega_j} (S, \R)$.  By Lemma \ref{lm:KZnorm2}, the estimate \pref{eq:recbound1} implies that
\begin{equation}
\label{eq:recbound2}
 \Vert   [z_j] \Vert_{\omega_j}  \leq     \Vert  [z_{j+1}] \Vert_{\omega_{j+1}} 
 \exp\bigl( \int_{t_j}^{t_{j+1}} \Lambda(g_s\omega) \,ds\bigr)  \, +\,   K^{(1)}_C \, \delta \,e^{2\tau_j} \,.
\end{equation}
For each $k\in \Z^+$, it follows by (reverse) induction on $1\leq j < k$ that
\begin{equation}
\label{eq:Hbound}
 \Vert   [z_j] \Vert_{\omega_j} \,\, \leq \,\,  2 K^{(1)}_C\, \delta^+  \, \Vert  [z_k] \Vert^+_{\omega_k} \,
  \exp\bigl( \int_{t_j}^{t_k} \Lambda(g_s\omega) \,ds\bigr)  \sum_{i=j}^{k-1}  \,e^{2\tau_i} \,.
\end{equation}
By the definition of the Hodge norm, for each $j\in \Z^+$, there exists a harmonic form $h_j \in
\Cal Z^{-1}_{\omega_j}(S)$ such that 
\begin{equation}
\label{eq:hrepr}
e_j = z_j -h_j  \in \Cal E^{-1}_{\omega_j} (S)  \quad \text{ and }   \quad 
\vert h_j \vert_{-1,\omega_j} \leq \Vert   [z_j] \Vert_{\omega_j} \,.
\end{equation}
 For each $j\in \Z^+$, let us define
 \begin{equation}
 \label{eq:Eid1}
 f_j= e_j - {\dc}_{-\tau_j} (e_{j+1})  \,\, \in \,\,\Cal E^{-1}_{\omega_j}(S)\,.
 \end{equation}
By the recursive identity \pref{eq:Eid1}  the following bound holds with respect to the Lyapunov norm $\Vert   \cdot \Vert_{-1}$ on the bundle of exact currents:
\begin{equation}
\label{eq:recbound3}
\Vert e_j \Vert_{-1, \omega_j} \leq  \Vert e_{j+1} \Vert_{-1, \omega_{j+1}}  +  
\Vert f_j \Vert_{-1, \omega_j}\,
\end{equation}
In fact, the restriction of the distributional cocycle $\{{\dc}_t\vert t\in\R\}$ to the bundle of exact currents is isometric with respect to $\Vert   \cdot \Vert_{-1}$.  For each $k\in \Z^+$, we derive from 
\pref{eq:recbound3}  by (reverse) induction on $1\leq j < k$ that
\begin{equation}
\label{eq:Ebound}
\Vert e_1 \Vert_{-1, \omega_1} \leq  \Vert e_k \Vert_{-1, \omega_k} + \sum_{j=1}^{k-1}
\Vert f_j \Vert_{-1, \omega_j} \,.  
\end{equation}
 By the splitting \pref{eq:gammasplit} and by the identities \pref{eq:hrepr} and \pref{eq:Eid1},
 it follows that
 \begin{equation}
  \label{eq:Eid2}
 f_j=  {\dc}_{-\tau_j} (h_{j+1}+ r_{j+1}) - (h_j +r_j)\,,
\end{equation}
hence by Lemma \ref{lm:Lnormineq}, by the bounds \pref{eq:cocycleestimates}, \pref{eq:Gammadeltabound}, \pref{eq:Hbound} and  \pref{eq:hrepr}, there exists a constant $K^{(2)}_C>1$ 
such that
\begin{equation}
\label{eq:remainderbound2}
\sum_{j=1}^{k-1} \Vert f_j \Vert_{-1, \omega_j} \leq K^{(2)}_C\, \delta^+ \, 
 \Vert  [z_k] \Vert^+_{\omega_k} \, \exp\bigl( \int_{t_0}^{t_k} \Lambda(g_s\omega) \,ds\bigr)  
 [ \sum_{j=1}^{k-1}  \,e^{2\tau_j} ]^2
\end{equation}
By the splitting \pref{eq:gammasplit} and by \pref{eq:Gammadeltabound}, \pref{eq:Hbound}, \pref{eq:hrepr}, \pref{eq:Ebound} and \pref{eq:remainderbound2}, there exists a constant $K^{(3)}_C >1$
such that for all $k>1$,
\begin{equation}
\label{eq:finalest}
\vert \gamma_1 \vert_{-1, \omega_1} \leq K^{(3)}_C\, \delta^+ \, 
 \vert \gamma_k \vert^+_{-1,\omega_k} \, \exp\bigl( \int_{t_1}^{t_k} \Lambda(g_s\omega) \,ds\bigr)  
 [ \sum_{j=1}^{k-1}  \,e^{2\tau_j} ]^2\,.
 \end{equation}
Finally, by the bound \pref{eq:cocycleestimates}, since $t_0=0$, 
\begin{equation}
\vert \gamma \vert_{-1, \omega} \leq e^{2(t_1-t_0)} \, \vert \gamma_1 \vert_{-1, \omega_1}  \,. 
\end{equation}
\end{proof}

\section{Proofs of main theorems}\label{proofs}

\subsection{Decomposing trajectories}
Let $\omega \in \Ab$ be an abelian differential and let $\{ \varphi_{\omega,t} 
\vert t\in \R\}$ denote the associated vertical translation flow, that is, the flow generated by a parallel
vertical normalized vector field with respect to the flat metric $R_\omega$ (with conical singularities) induced by $\omega$ on $S$. A point $x\in S$ will be called \emph{vertically non-singular} in the future [in the past] if the vertical forward orbit $\{\varphi_{\omega,t}(x) \vert t \geq 0\}$ [the vertical backward orbit $\{ \varphi_{\omega,t}(x) \vert t \leq 0 \}$] is well-defined (that is, it does not meet a singularity in finite time). Let us denote by $M_\omega^+$ [$M_\omega^-$] the full measure subset of points which are vertically non-singular in the future [in the past]. In the following we will consider the case of points vertically non-singular \emph{in the future}, since the case of points vertically non-singular in the past can be studied by similar arguments or it can be reduced to the previous case after a rotation of the abelian differential. Let $L_\omega$ denote the \emph{length functional }of the metric $R_\omega$. For any $(x,T)\in  M^+_\omega \times \R^+$, the path 
$$
\gamma_{\omega, x}(T)= \{ \varphi_{\omega,t}(x) \vert 0\leq t \leq T\}
$$ 
is a vertical non-singular geodesic for $R_\omega$ of length $L_\omega\bigr(\gamma_{\omega, x}(T)\bigl)=T$ starting at $x \in S$ and ending at $\varphi_{\omega,T}(x) \in S$.

 \smallskip
 By the Sobolev trace theorem (see for instance \cite{Adams}, Th. 5.4 (5)), any path $\gamma_{\omega, x}(T)$ induces by integration a $1$-dimensional current in $W^{-1}_\omega(S)$. Estimates on the deviation of ergodic averages of  functions in $H^{-1}_\omega(S)$ will be derived from the spectral gap estimates for the distributional cocycle $\{\dc_t \vert t\in \R\}$ proved in \S \ref{spectral}. In order to carry out this plan, we need to decompose any long non-singular vertical orbit segment $\gamma_{\omega, x}(T)$ into several sub-segments whose number and lengths are controlled in terms of the return times of the Teichm\"uller flow to a compact subset of the moduli space. The Sobolev norm of each of the subsegments will then be estimated by Lemma \ref{lm:iterativest} for an appropriate sequence 
 of return times which will be constructed below in \S  \ref{ss:sampling}. The following lemma allows 
 us to perform the above mentioned decomposition.

\begin{lemma}
\label{lm:chop} 
Let $\omega \in \Ab$ and let $\{T_k\}_{k \in \Z^+}$ be any non-decreasing divergent sequence of positive real numbers. For any $(x,T) \in M^+_\omega\times \R^+$, the vertical segment $\gamma_{\omega, x}(T)$ has a decomposition into consecutive sub-segments,
\begin{equation}
\label{eq:chop} 
\gamma_{\omega, x}(T) = \sum_{k=1}^n \sum_{m=1}^{m_k} \gamma_{\omega, x_{k,m}}(T_k)\,\, + \,\, 
\gamma_{\omega,y}(\tau) \,, 
\end{equation}
(with the convention that empty sums are equal to zero) such that  the numbers $m_1, \dots, m_k \in\Z^+\cup\{0\}$ and $\tau \geq 0$  satisfy the estimates
\begin{equation}
\label{eq:chopest}
m_k < T_{k+1} T_k^{-1} \quad \text{\rm and} \quad \tau < T_1\,.
\end{equation}
\end{lemma}
\begin{proof}  Let $n=\max\{k\in \Z^+ \,\vert \, T_k\leq T  \}$. The maximum exists (finite) since the sequence $\{T_k\}_{k \in \Z^+}$ is divergent. Let (for convenience) $T_0=0$.

We will construct a partition of the interval $[0,T]$ determined by  a finite set of appropriate times
 \begin{equation}
 \label{eq:partition}
\{ T_{k,m}\,\vert \, 1\leq k \leq n, 1\leq m\leq m_k\} \cup \{ T_y  \} \subset   [0, T]\,,
\end{equation} 
then define the required points on the orbit segment $\gamma_{\omega, x}(T)$ as
\begin{equation}
x_{k,m} = \varphi_{\omega, T_{k,m}} (x) \quad \text{\rm and } \quad y= \varphi_{\omega, T_y} (x)\,.
\end{equation}

Let $T_{n,1}= 0$ and let $m_n:= [ T T_n^{-1}] < T_{n+1}T_n^{-1}$. Define
$$
T_{n,m} =  (m -1)T_n  \,, \quad \text{ \rm for all }\, m\in \{1, \dots, m_n\}\,.
$$

The sequence $(T_{k,m})$ can then be constructed by a finite iteration of the 
following basic step. If $m_n, \dots, m_k \in \Z^+ \cup \{0\}$ with $m_k\geq 1$ and the set 
$$
\{T_{h,m} \vert  k\leq h \leq n\,, \,\, 1 \leq m\leq m_h \}
$$ 
has been constructed,  let 
\begin{equation}
\label{eq:defj}
j:= \max\{i\in \Z^+ \cup \{0\}  \vert  T_i \leq S_k:=T- (T_{k,m_k} + T_k)\}\,.
\end{equation}
Let  $m_h =0$ (that is, the set $\{T_{h,m}\}$ is empty) for all $j<h<k$ and let $m_j:= [S_kT_j^{-1}]\geq 1$. The set  $\{T_{j,1}, \dots, T_{j,m_j}\}$ is determined as follows:
\begin{equation}
\begin{aligned}
T_{j,1}&:=  T_{k,m_k} + T_k\,, \\ T_{j,m}&:= T_{j,1} + (m-1)T_j 
\,, \quad \text{ \rm for all }\, m\in \{1, \dots, m_j\}\,.
\end{aligned}
\end{equation}
It follows by \pref{eq:defj} that $T_{j+1} > S_k$, hence 
$$
m_j= [S_kT_j^{-1}] < T_{j+1}T_j^{-1}\,.
$$
Let $s:= \min\{k \in \Z^+\vert m_k \geq 1\}$ and let  $T_y= T_{s,m_s} + T_s$. By the definition of $m_s\geq 1$, we have $T_y + T_1 > T$, hence  $\tau:= T-T_y < T_1$. 

Let $k\in \{1, \dots, n\}$ be such that $m_k\geq 1$ and let $j\in\{1,\dots, k-1\}$ be the largest integer such that $m_j\geq 1$. By construction,  
\begin{equation}
\begin{aligned}
T_{k,m} + T_k &=  T_{k,m+1}\,, \quad &\text{ \rm  if }\,\, 1\leq m<m_k \,,  \\
T_{k,m} + T_k &=  T_{j,1}\,,\quad &\text{ \rm  if }\,\, 1\leq m=m_k \,,
\end{aligned}
\end{equation}
 It follows that $\{[T_{k,m}, T_{k,m}+T_k)  \, \vert \,1\leq k \leq n$, $1\leq m \leq m_k\}$ is a partition of
the interval $[0,T_y)$ and $[0,T_y) \cup [T_y, T_y +\tau]= [0,T]$. Finally, by the definition given above
of the points $x_{k,m}$ and $y$ belong to $\gamma_{\omega,x}(T) $, the decomposition \pref{eq:chop} 
holds, hence the argument is concluded.
\end{proof}
\begin{remark} In the application of Lemma \ref{lm:chop} (see below the proof of Theorem 
\ref{main}), we will consider a sequence $\{T_n\}$ of the form $T_n := e^{s_n}$ for an appropriate sequence $\{s_n\}$ of return times of the Teichm\"uller geodesic flow to a compact set in moduli space.
In this case, for all $k\in \{1,\dots,n\}$, 
$$
m_k \leq  e^{s_{k+1} -s_k}  \quad \text{ \rm and } \quad   \tau < e^{s_1}\,.
$$
\end{remark}

\subsection{A Sampling Lemma}
\label{ss:sampling}

The following `sampling lemma'  is motivated by the bounds proved in the spectral gap lemma for the distributional cocycle (Lemma \ref{lm:iterativest}) and  in the decomposition lemma for (vertical) trajectories of an abelian differential (Lemma \ref{lm:chop}), as it will be clear below in the proof of our main theorem.

\begin{lemma}\label{lm:sample} There exists a compact set $C \subset \Omega$ such that, for all
$\omega \in \Omega$ and almost all $\theta \in \Sone$, there is a diverging sequence $\{s_n : = s_n(\theta) \vert n\in \Z^+\}$ of forward times such that $g_{s_n} r_{\theta} \omega \in C$ and the following conditions hold:
\begin{enumerate}
\item for all $\epsilon >0$ there exists $n_\epsilon\in \Z^+$ such that for  $n\geq  n_\epsilon$, 
\begin{equation*}
0 < (s_{n+1} - s_n) \le \epsilon s_n \,;
\end{equation*}
\item for all $\lambda >0$ there is a $K_\lambda>0$ such that, for all $n\in \Z^+$,
\begin{equation*}
\sum_{k=1}^n e^{\lambda s_k} \le K_\lambda e^{\lambda s_n} \,.
\end{equation*}
\end{enumerate}
\end{lemma}
\begin{proof} We will construct this sequence following the paper~\cite{Athreya}. All theorem and equation numbers in the rest of the paragraph are drawn from~\cite{Athreya}. Let $V: \Omega \rightarrow \R$ be as in the statement of Lemma 2.10 (to be precise, fix some $\delta >0$, and let $V : = V_{\delta}$). Let $C_l : = \{\omega \in \Omega: V(\omega) \le l\}$. The family $\{C_l\}_{l>0}$ forms an exhaustion of $\Omega$ by compact sets. The function $V$ is $SO(2,\R)$-invariant, thus, $r_{\theta} \omega \in C_l$ iff $\omega \in C_l$. As a consequence, we can work on the quotient $SO(2,\R)\backslash SL(2, \R) \omega$, which we can view as the hyperbolic plane $\h^2$.  For any $\omega\in\Omega$, we identify the $SO(2, \R)$-orbit $\{r_{\theta}\omega \vert \theta \in \Sone\}$ with the point $i \in \h^2$, and $\omega$ with the upward pointing unit tangent vector. The point $g_t r_{\theta} \omega$ can be identified with $i \cdot g_t r_{\theta}$. 

Let $d$ denote the Gromov constant for $\h^2$ (the smallest $c$ such that any side of a triangle is contained in a $c$-neighborhood of the other two). Let $l_0$ be defined as in the statement of Theorem 2.1, and $l>l_0$ so that 
\begin{equation*}
2d < d(C_l^c , C_{l_0}) := \inf \{ d(z_1, z_2) \,\vert \,z_1,\, z_2 \in \h^2, \,z_1 \notin C_l 
\text{ and } z_2 \in C_{l_0} \}\,.
\end{equation*}
(The points $z_1$ and $z_2 \in \h^2$ represent $SO(2, \R)$-orbits within the $SL(2, \R)$-orbit
of a given $\omega \in \Omega$).

Let $0<s <d(C_l^c , C_{l_0})-2d$, let $s_0 (\theta) = 0$ and define
\begin{equation*}
s_n(\theta) =  s_{n-1}(\theta) + s
\end{equation*}
 if $i\cdot g_t r_{\phi} \in C_l$, for all $t \in (s_{n-1}(\theta), s_{n-1}(\theta) + s]$ and for all 
 $\phi \in \Sone$  such that $d(i\cdot g_t r_{\phi}, i\cdot g_t r_{\theta}) < d$, and otherwise set
\begin{displaymath} 
s_n(\theta) =  \inf \{t>s_{n-1}(\theta)+s\,\vert \,\exists \phi \mbox{ s. t. }d(i\cdot g_t r_{\phi}, i \cdot g_t
r_{\theta}) < d, i\cdot g_t r_{\phi} \in C_{l_0}\} .
\end{displaymath} 
That is, if all nearby trajectories remain in $C_l$ for the interval of time of length $s>0$ after our previous sampling time $s_{n-1}$, we sample at $s_{n-1}+s$. If not, we sample at the first time after $s_{n-1}+s$ when a nearby trajectory re-enters $C_{l_0}$. This scheme guarantees that $g_{s_n} r_{\theta} \omega \in C_l$ and that $s_{n+1} \geq s_n +s$, for all $n\in \N$, hence condition $(2)$ is automatically satisfied. In fact, since for each $j\in \N$, the cardinality of $\{k\in\N \,\vert\, [s_k]=j\}$ is at most 
$1+[1/s]$,
\begin{equation*}
\sum_{k=1}^n e^{\lambda s_k} \leq \sum_{k=1}^n  e^{\lambda ([s_k] +1)} \leq   
(1+[\frac{1}{s}])\,\sum_{j=1}^{[s_n]+1}  e^{\lambda j}  \leq K_{\lambda,s}\, e^{\lambda s_n} \,.
\end{equation*}

It remains to check condition $(1)$. In~\cite{Athreya}, the following lemma was proved

\begin{lemma} Fix notation as above. Let $t_0(\theta)=0$ and, for all $k >0$, let
\begin{displaymath} 
t_{2k}(\theta) = \inf \{t>t_{2k-1}\,\vert\, \exists \phi \mbox{ such that }d(i\cdot g_t r_{\phi}, i\cdot g_t
r_{\theta}) < d, i\cdot g_t r_{\phi} \in C_{l_0}\}
\end{displaymath} 
and
\begin{displaymath} t_{2k+1}(\theta) = \inf \{t> t_{2k}\,\vert\, \exists \phi
\mbox{ such that } d (i\cdot g_t r_{\phi}, i\cdot g_t r_{\theta}) < d,
i\cdot g_t r_{\phi} \notin C_l\}.
\end{displaymath} 
Let  $\tau_k : = t_k - t_{k-1}$. The Lebesgue measure of the set $\{\theta\in \Sone \,\vert \,\tau_{2k}(\theta) > t \}$ decays exponentially in $t\in \R^+$ (uniformly with respect to the basepoint  $\omega\in \Omega$ and $k\in \N$). 
\end{lemma}

This result is implicitly contained in the proof of Theorem 2.3 (in~\cite{Athreya}). We would like to thank an anonymous referee for pointing out a small mistake in that proof, which we indicate how to correct here. It was wrongly assumed there that the times $t_k(\theta)$ are locally constant in $\theta\in \Sone$ in order to construct ``intervals'' $I_k(\theta) 
= \bigcap_{i=1}^{k}(t_i^{\prime})^{-1} (t_i^{\prime}(\theta))$, where the $t_i^{\prime}$'s are an auxiliary set of times (for the precise definition, see~\cite{Athreya}, page 139). To fix the argument, we replace the $I_k(\theta)$'s (which are not really intervals) with $I_k^{\prime}(\theta)$, defined as follows. Fix $\epsilon >0$, let 
$$
I_k^{\prime}(\theta) = \bigcap_{i=1}^{k}(t_i^{\prime})^{-1}(t_i(\theta) 
-\epsilon/2^i, t_i(\theta) +\epsilon/2^i) \,.
$$
The rest of the proof in~\cite{Athreya} then goes through with minor modifications. We now continue with the proof of Lemma~\ref{lm:sample}.

We claim that for  each $n\in\N$ such that $s_{n+1} \neq s_n + s$, there exists $k\in\N$ such that 
$t_{2k-1} < s_n +s$ and $s_{n+1} = t_{2k}$. In fact, let $k\in \N$ be the largest integer such that
$t_{2k-2} \leq s_n$. Since $s_{n+1} \neq s_n + s$ and $2d+s < d(C_l^c , C_{l_0})$, by the above construction  $t_{2k-1} < s_n +s$ and $t_{2k} > s_n +s$. Thus by definition $s_{n+1} = t_{2k}$ as claimed. 

Let $B_n(\epsilon)$ be the set  of $\theta \in \Sone$ such that  $s_{n+1}(\theta) - s_n(\theta)
 \geq \epsilon s_n(\theta)$ and let $B(\epsilon)$ be the set of $\theta \in \Sone$  such that 
 $\theta \in B_n(\epsilon)$ for infinitely many $n\in \N$.  By the above claim, it follows that 
 $s_{n+1} - s_n \geq \epsilon s_n$ implies either that $s\geq \epsilon s_n$ or that there exists
 $k\in \N$ such that $\tau_{2k} \geq \epsilon s_n -s$. Since $s_n \geq s n$, the first case holds for 
 (at most) finitely many $n\in \N$, while in the second case the Lebesgue measure of $B_n(\epsilon)$ 
 is exponentially small in $s n-s$. Thus the series of the Lebesgue measures of the set $B_n(\epsilon)$ 
 is convergent and, by the Borel-Cantelli law, the set $B(\epsilon)$ has measure $0$.

Let $B = \cup_{j \in \N} B(1/j)$. The set $B$ has measure $0$ and condition $(1)$ holds on the
complement $\Sone\setminus B$.
 \end{proof}
 
 \subsection{Conclusion}
 Let $\omega \in \Ab$. By the Sobolev trace theorem \cite{Adams}, Th. 5.4 (5), any smooth path in $S$
 induces by integration a ($1$-dimensional) current which belongs to $W^{-1}_\omega(S)$.
 For any geodesic segment on the flat surface $(S, R_\omega)$ there exists a precise bound,
 proved in  \cite{Forni}, of the Sobolev norm of the induced current in terms of the length of the 
 segment and of the geometry of the flat surface. 
 
 Let $L_\omega$ denote the \emph{length functional} of the flat metric $R_\omega$ and let $\vert\!\vert\omega \vert\!\vert$ denote the $R_\omega$-length of the shortest \emph{saddle connection} (=geodesic connecting two zeros with no other zeros in the interior).  

\begin{lemma}
\label{lm:Sobtrace} (Lemma 9.2 in \cite{Forni}) There exists a constant $K_\Omega>1$ such that, for all $\omega \in \Omega$ and for any non-singular geodesic segment $\gamma \subset (S,R_\omega)$, 
\begin{equation}
\label{eq:Sobtrace}
\vert \gamma \vert_{-1,\omega} \,\leq\, 
K_\Omega(1+ {\frac{L_\omega(\gamma)}{\vert\!\vert\omega\vert\!\vert}})\,.
\end{equation}
\end{lemma} 
The above Lemma was proved in \cite{Forni} for horizontal or vertical segments but it can be
extended by rotation invariance to any non-singular geodesic segment. In fact, $L_{r_\theta\omega} = L_\omega$ and $\vert\!\vert r_{\theta}\omega\vert\!\vert = \vert\!\vert\omega \vert\!\vert$ for all $\theta \in \Sone$ since the flat metric $R_\omega$ is invariant under the action of $SO(2,\R)$ on $\omega\in \Omega$. 

For any $\omega\in \Omega$, let $d_\omega$ denote the diameter of the flat surface $M_\omega$.
 For any (piecewise geodesic) path $\gamma$ in $(S,R_\omega)$, there exists a \emph{closed }(piecewise geodesic) path $\bar{\gamma}$ such that 
 $L_\omega(\bar{\gamma}-\gamma) \leq d_\omega$. The closed path 
 $\bar{\gamma}$ is obtained by joining the endpoints of $\gamma$ by a geodesic path of minimal
length in $(S,R_\omega)$.  For any compact set $C\subset \Omega$ there exist $d_C$, $s_C >0$ such that $d_\omega\leq d_C$ and $\vert\!\vert\omega\vert\!\vert \geq s_C$ for all $\omega\in C$. Hence by Lemma \ref{lm:Sobtrace} there exists $\delta_C>0$ such, that for all $\omega\in C$,
\begin{equation}
\label{eq:deltabound}
\vert \bar{\gamma}-\gamma \vert_{-1, \omega} \leq  K_\Omega(1+ 
{\frac{d_\omega}{\vert\!\vert\omega\vert\!\vert}}) \leq \delta_C\,.
\end{equation}
It follows that any (geodesic) path $\gamma \in \Gamma_C(\delta)$ for any $\delta \geq \delta_C$.
In fact, for all $t\in \R$, the current ${\dc}_t(\gamma)$ is represented by a (geodesic) path and
 any closed (geodesic) path belongs to $\Cal Z^{-1}_\omega(S)$. Hence, if $g_t\omega \in C$,
by the bound \pref{eq:deltabound} the distance $\delta_g({\dc}_t(\gamma))$ from the sub-bundle of closed currents is $\leq \delta_C$.

Smooth $1$-forms on naturally define $1$-currents. In fact, there is a standard embedding of the space of $1$-forms into  the space of $1$-currents induced by the following bilinear pairing: for any 
pair $(\xi_1, \xi_2)$ of $1$-forms on $S$, 
 $$
 \<\xi_1, \xi_2\> = \int_S \xi_1 \wedge \xi_2 \,.  
 $$
 It follows from the definitions that the forms $\operatorname{Re} (\omega) $ and $\operatorname{Im} (\omega)$ belong as distributions to the Sobolev space $W^{-1}_\omega(S)$, for any  $\omega\in \Ab$.

\begin{proof} [Proof of Theorem \ref{main}] Let $\omega \in \Omega$. We claim that
there exist a real number $\alpha:= \alpha(\Omega)>0$ and a measurable function 
$K: \Sone \to \R^+ \cup \{+\infty\}$, almost everywhere finite, such that for all $\theta \in \Sone$,
for all $x\in S$ vertically non-singular for the abelian differential $r_\theta \omega\in \Omega$ 
and all $T\geq 1$,
\begin{equation}
\label{eq:mainclaim}
\vert \gamma_{r_\theta\omega, x}(T) - T \, \operatorname{Re} (r_\theta\omega) \vert_{-1, \omega} 
\leq   K(\theta) \, T ^{1-\alpha} \,.
\end{equation}
For almost all $\theta \in \Sone$, let $s_n:=s_n(\theta)$ be the sequence of visiting
times of the orbit $\{g_t r_\theta \omega \,\vert \, t\geq 0\}$ to a compact set $C\subset
\Omega$ constructed in Lemma \ref{lm:sample} and let $\{T_n\}$ be the sequence defined 
as $T_n:= e^{s_n}$, for all $n\in \Z^+$.

Let us consider the decomposition \pref{eq:chop}, determined by the sequence $\{T_n\}$ as in Lemma \ref{lm:chop}, of the orbit segment $\gamma_{r_\theta\omega, x}(T)$ into consecutive sub-segments 
and let us introduce the currents:
\begin{equation}
\begin{aligned}
\gamma_{k,m}&= \gamma_{r_\theta\omega,x_{k,m}}(T_k)  - T_k  \, \operatorname{Re} (r_\theta\omega) \,, \\
\gamma_0&= \gamma_{r_\theta\omega,y}(\tau) - \tau  \, \operatorname{Re} (r_\theta\omega) \,.
\end{aligned}
\end{equation}
It follows that there is a decomposition
\begin{equation}
\label{eq:finalchop}
\gamma_{r_\theta\omega,x}(T) - T \, \operatorname{Re} (r_\theta\omega) = \sum_{k=1}^n \sum_{m=1}^{m_k}
\gamma_{k,m}\, + \, \gamma_0  \,\, \in \,\, W^{-1}_\omega(S)\,.
\end{equation}
Since $\tau \leq e^{s_1}$ by Lemma \ref{lm:chop}, it follows from Lemma \ref{lm:Sobtrace} that
\begin{equation}
\label{eq:rest}
\vert \gamma_0 \vert_{-1, \omega} \leq  2 K_\Omega
(1+ {\frac{e^{s_1}}{\vert\!\vert\omega\vert\!\vert}})\,.
\end{equation}
The Sobolev norms of each current $\gamma_{k,m} \in W^{-1}_\omega(S)$ can be estimated by
applying Lemma \ref{lm:iterativest}. In fact, by definition $\gamma_{k,m} \in \Cal I^{\perp}_\omega(S)$,
that is,
\begin{equation}
\begin{aligned}
\<\gamma_{k,m} \wedge \operatorname{Re}(\omega),1\>&= \<\gamma_{k,m} , \operatorname{Re}(\omega)\> =0  \,; \\
\< \gamma_{k,m} \wedge \operatorname{Im}(\omega),1\> &= \<\gamma_{k,m} , \operatorname{Im}(\omega)\> =0  \,.
\end{aligned}
\end{equation}
Since $\operatorname{Re} (r_\theta\omega) \in \Cal Z^{-1}_\omega(S)$, by the above discussion, in particular by formula \pref{eq:deltabound}, the currents $\gamma_{k,m} \in \Gamma_C(\delta)$ for any $\delta\geq \delta_C$. 

For all $n\in \Z^+$, let $\omega_n= g_{s_n} r_\theta \omega$. Let us verify that there exists a constant $K_1=K_1(C)>0$ such that, for all $1\leq k \leq n$ and $1\leq m\leq m_k$, 
\begin{equation}
\label{eq:dcskbound}
\vert \dc_{s_k} ( \gamma_{k,m} ) \vert_{-1, \omega_k} \leq K_1 \,.
\end{equation}
In fact, since by construction $\gamma_{r_\theta\omega,x_{k,m}}(T_k) $ is a vertical segment of length $T_k=e^{s_k}$ for the differential $r_\theta\omega\in \Omega$, it is a vertical segment of length exactly equal to $1$ for $\omega_k= g_{s_k} r_\theta \omega$, hence its Sobolev norm is bounded by Lemma 
\ref{lm:Sobtrace}. By the definition of the distributional cocycle $\{\dc_t\vert t\geq 0\}$,  it follows
that there exists a constant $K_0= K_0(C)>0$ such that
$$
\vert {\dc}_{s_k} \left(\gamma_{r_\theta\omega,x_{k,m}}(T_k) \right) \vert_{-1, \omega_k}
\,\leq \, K_0\,.
$$
On the other hand, since ${\dc}_{s_k}\left( \operatorname{Re} (r_\theta\omega) \right) = e^{-s_k} \, 
\operatorname{Re} (\omega_k)$, it follows that 
$$
\vert {\dc}_{s_k} \left( T_k \operatorname{Re} (r_\theta\omega) \right)  \vert_{-1, \omega_k}  = 1 \,.
$$
By the two above bounds we obtain inequality \pref{eq:dcskbound} with $K_1=K_0+1$.

Thus by Lemma \ref{lm:recurrest} and Lemma \ref{lm:iterativest}, for any $\lambda<\lambda_1<1$ 
there exists a measurable function $K_2: \Sone \to \R^+$ such that for almost all $\theta \in \Sone$, 
\begin{equation}
\vert \gamma_{k,m}  \vert _{-1, \omega} \leq K_2(\theta) \, \delta^+ \, e^{\lambda_1 s_k} 
 [ \sum_{j=0}^{k-1}  e^{2(s_{j+1}-s_j)} ]^3 \,.
\end{equation}
By conditions $(1)$ and $(2)$ in Lemma \ref{lm:sample}, for every $\epsilon>0$
there exists a measurable function $K_3: \Sone \to \R^+$ such that, for almost all $\theta \in \Sone$, 

$$
\sum_{j=0}^{k-1}  e^{2(s_{j+1}-s_j)}  \leq K_3(\theta) \,e^{\epsilon s_k} \,, \quad \text{ \rm for all }
\, k\in\Z^+  \,.
$$
Hence, for any $\lambda_1 < \lambda_2<1$ there is a measurable function $K_4: \Sone \to \R^+$ such that, for almost all $\theta \in \Sone$, 
\begin{equation}
\label{eq:piecesest}
\vert \gamma_{k,m}  \vert _{-1, \omega} \leq K_4(\theta) \, \delta^+ \, e^{\lambda_2 s_k}  \,.
\end{equation}
Finally, by the decomposition \pref{eq:finalchop} and the bounds \pref{eq:rest} and 
\pref{eq:piecesest}, there exists  a measurable function $K_5: \Sone \to \R^+$ such that, for 
almost all $\theta \in \Sone$,  
\begin{equation}
\label{eq:gatherest1}
\vert \gamma_{r_\theta\omega, x}(T) - T \, \operatorname{Re} (r_\theta\omega)   \vert _{-1, \omega} 
\leq K_5(\theta)  \sum_{k=1}^n  m_k \,  e^{\lambda_2 s_k}\,.
\end{equation}
Since  $m_k \leq e^{s_{k+1} -s_k}$, for all $1\leq k \leq n$, and $e^{s_n} =T_n \leq T$ in the decomposition \pref{eq:finalchop}, by Lemma \ref{lm:chop}, it follows  from conditions $(1)$ and $(2)$ in the sampling Lemma \ref{lm:sample} that, for all  $\lambda_2< \lambda_3 <1$, there exists a  measurable function $K_6: \Sone \to \R^+$ 
such that, for almost all $\theta \in \Sone$,
\begin{equation}
\label{eq:gatherest2}
\sum_{k=1}^n  m_k \,  e^{\lambda_2 s_k}   \leq    K_6(\theta)  \,  e^{\lambda_3 s_n}  
 \leq   K_6(\theta)  \,  T^{\lambda_3} \,.
\end{equation}
The proof of the claim \pref{eq:mainclaim} follows from \pref{eq:gatherest1} and \pref{eq:gatherest2}.

Let $f \in H^1(S) = H^1_\omega(S)$ and let $\psi_f = f\,\operatorname{Im}(r_\theta\omega) \in W^1_\omega(S)$.
Since
$$
\int_0^T f \left (\varphi_{r_\theta\omega, s} (x)\right) ds  -  T \int_S f\, dA_\omega 
= \<\gamma_{r_\theta\omega, x}(T) - T \, \operatorname{Re} (r_\theta\omega) , \psi_f \>\,, 
$$
it follows from the estimate \pref{eq:mainclaim} that 
$$
\left \vert \int_0^T f \left (\varphi_{r_\theta\omega, s} (x)\right) ds  -  T \int_S f\, dA_\omega \right\vert
\leq K(\theta) \, \vert f \vert_{1, \omega} \, T^{1-\alpha}\,.
 $$
 \end{proof}
\begin{proof} [Proof of Theorem \ref{mainhom}]  Let $h_\theta\in H_1(S,\R)$
be the Poincar\'e dual of the cohomology class  $[\operatorname{Re}(r_\theta\omega)] \in H^1(S, \R)$. 
Any cohomology class $c\in H^1(S, \R)$ can be represented as $c= [ \psi_c]$ for a smooth 
\emph{closed} $1$-form $\psi_c$ supported on $S\setminus \Sigma_\omega$, which therefore 
belongs to $W^1_{\omega}(S)$. The following identity holds for any $c\in H^1(S,\R)$:
$$
\<c, h_{r_\theta\omega,x}(T) - T h_\theta\> = \< \bar{\gamma}_{r_\theta\omega,x}(T) 
-T \operatorname{Re}(r_\theta\omega), \psi_c \>\,.
$$
Hence Theorem \ref{mainhom} also follows from \pref{eq:mainclaim}.

\end{proof}

\noindent\textbf{Acknowledgements:} We would like to thank the anonymous referees for their careful study of and valuable comments on the paper. We are also grateful to E. Gutkin for pointing out to us
an inaccuracy in the presentation of rational polygonal billiards in \S \ref{RPB}.

\end{document}